\documentclass{gtart_h}
\def\ifplaintex{\expandafter\ifx\csname documentclass\endcsname\relax}

\def\gtp{{\mathsurround=0pt\it $\cal G\mskip-2mu$eometry \&\ 
$\cal T\!\!$opology $\cal P\!$ublications}}  % GT publications

\def\Addressesr{\bigskip
{\small \parskip 0pt \leftskip 0pt \rightskip 0pt plus 1fil \def\\{\par}
\sl\theaddress\par
\medskip
\rm Email:\stdspace\tt\theemail\hfill\rm Received:\qua\receiveddate \par}}

\def\recd{{\small Received:\qua\receiveddate\ifx\reviseddate\relax
\else\qquad Revised:\qua\reviseddate\fi\par}} 

\def\lognumber#1{\def\thelognumber{#1}}
\def\volumenumber#1{\def\thevolumenumber{#1}}
\def\volumeyear#1{\def\thevolumeyear{#1}}
\def\papernumber#1{\def\thepapernumber{#1}}
\def\pagenumbers#1#2{\def\startpage{#1}\def\finishpage{#2}}
\def\published#1{\def\publishdate{#1}}

\def\received#1{\def\receiveddate{#1}}

\def\accepted#1{\def\accepteddate{#1}}
\def\asciititle#1{\def\theasciititle{#1}}

\def\asciikeywords#1{\def\theasciikeywords{#1}}

\let\\\par\let\thelognumber\relax\let\thevolumenumber\relax
\let\thepapernumber\relax\let\thevolumeyear\relax\let\startpage\relax
\let\finishpage\relax\let\publishdate\relax\let\receiveddate\relax
\let\reviseddate\relax\let\accepteddate\relax\let\theasciititle\relax
\let\theasciiauthors\relax
\let\theasciiabstract\relax\let\theasciikeywords\relax

\let\theasciiemail\relax

\ifplaintex
\font\logobig=cmssbx10 scaled 3836
\font\logomed=cmssbx10 scaled 2557
\else
\font\logobig=cmssbx10 scaled 4200
\font\logomed=cmssbx10 scaled 2800
\fi

\long\def\makeagttitle{   %%% start of definition of \makeagttitle
\count0=\startpage
\agt\hfill      %   Journal title (top left) 
\hbox to 45truept{\vbox to 0pt{\vglue -13truept{\logomed A\kern -.37em{\logobig 
T}\kern -.38em G}\vss}\hss}
\break
{\small Volume \thevolumenumber\ (\thevolumeyear)
\startpage--\finishpage\nl
Published: \publishdate}

\vglue .25truein

{\parskip=0pt\leftskip 0pt plus
1fil\def\\{\par\smallskip}{\Large\bf\thetitle}\par\medskip} \vglue
0.05truein

{\parskip=0pt\leftskip 0pt plus 1fil\def\\{\par}{\sc\theauthors}
\par\medskip}%
 
\vglue 0.03truein 

{\small\leftskip 25truept\rightskip 25truept{\bf Abstract}\stdspace\theabstract

{\bf AMS Classification}\stdspace\theprimaryclass
\ifx\thesecondaryclass\relax\else; \thesecondaryclass\fi\par
{\bf Keywords}\stdspace \thekeywords\par}\vglue 7truept

}   %%%% end of definition of \makeagttitle

\ifplaintex
\font\phead=cmsl9 scaled 950
\font\pnum=cmbx10 scaled 913
\font\pfoot=cmsl9 scaled 950
\headline{\vbox to 0pt{\vskip -4.5mm\line{\small\phead\ifnum
\count0=\startpage ISSN 1472-2739 (on-line) 1472-2747 (printed)
\hfill {\pnum\folio}\else\ifodd\count0\def\\{ }% 
\ifx\theshorttitle\relax\thetitle\else\theshorttitle\fi\hfill{\pnum\folio}
\else\def\\{ and }{\pnum\folio}\hfill\ifx\theshortauthors\relax\theauthors
\else\theshortauthors\fi\fi\fi}\vss}}
\footline{\vbox to 0pt{\vglue 0mm\line{\small\pfoot\ifnum\count0=\startpage
\copyright\ \gtp\hfill\else
\agt, Volume \thevolumenumber\ (\thevolumeyear)\hfill\fi}\vss}}
\else
\font=cmsl9 scaled 1050
\font\lnum=cmbx10 
\font=cmsl9 scaled 1050
\makeatletter
\def\@oddhead{{\small\lhead\ifnum\count0=\startpage ISSN 1472-2739 
(on-line) 1472-2747 (printed)\hfill {\lnum\number\count0}\else\ifodd\count0
\def\\{ }\ifx\theshorttitle\relax \thetitle \else\theshorttitle\fi\hfill
{\lnum\number\count0}\else\def\\{ and }{\lnum\number\count0}
\hfill\ifx\theshortauthors\relax 
\theauthors\else\theshortauthors\fi\fi\fi}}\def\@evenhead{\@oddhead}
\def\@oddfoot{\small\lfoot\ifnum\count0=\startpage\copyright\ \gtp\hfill\else
\agt, Volume \thevolumenumber\ (\thevolumeyear)\hfill\fi}
\def\@evenfoot{\@oddfoot}
\makeatother
\fi
\let\maketitlepage\makeagttitle
\let\makeshorttitle\maketitlepage
\let\maketitle\maketitlepage

\newwrite\gtoutfile
\long\gdef\makeheadfile{  %%% start of definition of \makeheadfile
{\def\\{, }\def\s{ }
\immediate\openout\gtoutfile head.xxx
\immediate\write\gtoutfile{Proxy-for: \ifx\theasciiauthors\relax
\theauthors\else\theasciiauthors\fi\s<\ifx\theasciiemail\relax\theemail\else\theasciiemail\fi>}
\immediate\write\gtoutfile{\noexpand\\}
\immediate\write\gtoutfile{Authors: \ifx\theasciiauthors\relax
\theauthors\else\theasciiauthors\fi}
{\def\\{ }\immediate\write\gtoutfile{Title: \ifx\theasciititle\relax
\thetitle\else\theasciititle\fi}}
\immediate\write\gtoutfile{Subj-class: GT or SG, GR etc}
\immediate\write\gtoutfile{MSC-class: \theprimaryclass\ifx\thesecondaryclass\relax\else, \thesecondaryclass\fi}
\immediate\write\gtoutfile{Journal-ref: Algebr. Geom. Topol. \thevolumenumber\s
(\thevolumeyear) \startpage-\finishpage}
\immediate\write\gtoutfile{Comments: Published by Algebraic and
Geometric Topology at}
\immediate\write\gtoutfile{\s\s\s  http://www.maths.warwick.ac.uk/agt/AGTVol\thevolumenumber/agt-\thevolumenumber-\thepapernumber.abs.html}
\immediate\write\gtoutfile{\noexpand\\}
\immediate\write\gtoutfile{}
\ifx\theasciiabstract\relax
\immediate\write\gtoutfile{\theabstract}\else
\immediate\write\gtoutfile{\theasciiabstract}\fi
\immediate\write\gtoutfile{}
\immediate\write\gtoutfile{\noexpand\\}
\immediate\write\gtoutfile{}
\immediate\closeout\gtoutfile}}  %%% end of definition of \makeheadfile

\def\maketitlepage{\makeagttitle\makeheadfile}
\let\makeshorttitle\maketitlepage
\let\maketitle\maketitlepage

\lognumber{59}
\volumenumber{5}
\volumeyear{2005}
\papernumber{59}
\pagenumbers{1481}{1503}
\received{19 November 2004} 
\accepted{20 October 2005}
\published{1 November 2005}

\usepackage{amsmath,amssymb,epic,eepic}

\def\figref#1{\hyperlink{#1anchor}{Figure~\ref*{#1}}}
\def\anchor#1{\noindent\hypertarget{#1anchor}{\smash{$\phantom{99}$}}}

\newcommand{\Z}{\ensuremath{\mathbb{Z}}}
\newcommand{\R}{\ensuremath{\mathbb{R}}}
\newcommand{\GT}{\ensuremath{\mathcal{X}}} % space of G-trees
\newcommand{\GM}{\ensuremath{\mathcal{M}}} % space of G-morphisms
\newcommand{\GE}{\ensuremath{\mathcal{E}}} % space of equivarant maps
\newcommand{\GS}{\ensuremath{\mathcal{T}}} % space of transverse maps
\newcommand{\GL}{\ensuremath{\mathcal{G}}} % space of measured laminations
\newcommand{\Do}{\ensuremath{\mathcal{D}o}} % domain function
\newcommand{\Ra}{\ensuremath{\mathcal{R}a}} % range function
\newcommand{\F}{\ensuremath{\mathcal{F}}} % partition
\newcommand{\GD}{\ensuremath{\mathcal{D}}} % deformation space
\newcommand{\GSs}{\ensuremath{\mathcal{S}}} % set of G-spaces 
\newcommand{\GMs}{\ensuremath{\mathcal{S'}}} % set of G-maps 
\newcommand{\CC}{\ensuremath{\mathcal{C}}} % conjugacy class

\newtheorem{theorem}{Theorem}[section]
\newtheorem{lemma}[theorem]{Lemma}
\newtheorem{proposition}[theorem]{Proposition}

\newtheorem{definition}[theorem]{Definition}

\newtheorem{remark}[theorem]{Remark}

\DeclareMathOperator{\Out}{Out}
\DeclareMathOperator{\gr}{graph}
\DeclareMathOperator{\vol}{vol}
\DeclareMathOperator{\length}{length}

\begin{document}

%%%%%%%%%%%%%%%%%%%%%%%%%%%%%%%%%%%%%%%%%%%%%%%%%%%%%%%%%%%%%%%%%%%%%%%%%%%%% 

\title{Contractibility of deformation spaces of $G$--trees}
\asciititle{Contractibility of deformation spaces of G-trees}
\author{Matt Clay} 
\address{Department of Mathematics, University of Utah\\Salt 
Lake City, UT 84112-0090, USA}
\email{clay@math.utah.edu}

\begin{abstract}
Forester has defined spaces of simplicial tree actions for a finitely
generated group, called deformation spaces.  Culler and Vogtmann's
Outer space is an example of a deformation space.  Using ideas from
Skora's proof of the contractibility of Outer space, we show that
under some mild hypotheses deformation spaces are contractible.
\end{abstract}

\asciikeywords{G-tree, deformation space, Outer space}

\maketitle

%%%%%%%%%%%%%%%%%%%%%%%%%%%%%%%%%%%%%%%%%%%%%%%%%%%%%%%%%%%%%%%%%%%%%%%%%%%%% 

Culler and Vogtmann's Outer space is a good geometric model for
$\Out(F_n)$, the outer automorphism group of a finitely generated free
group of rank $n\geq 2$, for three reasons:

\begin{itemize}

\item[(1)] Outer space is contractible;

\item[(2)] point stabilizers are finite; and

\item[(3)] there is a equivariant deformation retract on which the
action is cocompact \cite{ar:CV}.

\end{itemize}

\noindent
Outer space is the analog of Teichm\"uller space for the mapping class
group of a closed negatively curved surface or of the symmetric space
for an arithmetic group.  See \cite{ar:B} and \cite{ar:V} for a survey
of some results about $\Out(F_n)$ obtained from using this connection
between the three classes of groups.  Also see \cite{ar:BV} for some
open questions about the similarities and differences.

Recall that Outer space is the moduli space of free actions of a free
group on a simplicial tree.  Forester has defined a generalization of
Outer space for an arbitrary finitely generated group $G$
\cite{ar:F1}.  The generalization allows actions which are not free but
requires the subgroups with fixed points to be the same among all
actions in the moduli space.  Unfortunately these spaces are not
$\Out(G)$--invariant in general.  Nevertheless, in the cases when the
space is invariant under $\Out(G)$ these spaces have the potential to
provide information about the structure of $\Out(G)$.  The purpose of
this paper is to show that these spaces share the first of the above
mentioned properties with Outer space, i.e.\ they are contractible.

For a finitely generated group $G$, a $G$--\textit{tree} is a metric
simplicial tree on which $G$ acts by isometries.  Two $G$--trees $T$
and $T'$ are equivalent if there is a $G$--equivariant isometry
between then.  When we speak of a $G$--tree we will always mean the
equivalence class of the $G$--tree.  A subgroup is called an
\textit{elliptic} subgroup for $T$ if it has a fixed point in $T$.
Given a $G$--tree there are two moves one can perform to the tree that
do not change whether or not subgroups of $G$ are elliptic.  These
moves correspond to the isomorphism $A \cong A *_C C$ and are called
\textit{collapse} and \textit{expansion}.  For a detailed description
of the moves see \cite{ar:F1}.  In \cite{ar:F1} Forester proves the
converse, namely if two cocompact $G$--trees have the same elliptic
subgroups, then there is a finite sequence of collapses and expansions
(called an \textit{elementary deformation}) transforming one $G$--tree
to the other.  A $G$--tree $T$ is \textit{cocompact} if the quotient
$T/G$ is a finite graph.

We let \GT \ denote a maximal set of cocompact $G$--trees which are
related by an elementary deformation.  By the theorem of Forester
mentioned above, an equivalent definition is as the set of all
cocompact $G$--trees that have the same elliptic subgroups as some
fixed $G$--tree.  Both of these interpretations are utilized in the
following.  This set \GT \ is called a \textit{unnormalized
deformation space}.  We will always assume that the $G$--trees are
minimal, irreducible and that $G$ acts without inversions.  See
section \ref{sc:pre} for these definitions.

As is common practice in spaces of this nature, we projectivize by
taking the quotient of \GT \ under the action of $\R^+$ by homothety.
The quotient $\GT / \R^+$ is called a \textit{deformation space} and
is denoted \GD.  Outer space is an example of a deformation space for
a finitely generated free group where the only elliptic subgroup is
the trivial group.  Culler and Vogtmann described a contraction of the
spine of Outer space using combinatorial methods and a ``Morse-like''
function \cite{ar:CV}.  Skora showed in a different manner that Outer
space is contractible \cite{ar:S}.  The method of Skora is to homotope
the unnormalized deformation space projecting to Outer space to a set
homeomorphic to a simplex$\times \R^+$ by continuously unfolding
$G$--trees in the unnormalized deformation space.  This homotopy
descends to Outer space, proving its contractibility.  It is this idea
which we extend to show:

\medskip
\textbf{Theorem \ref{th:d_contractible}}\qua \textsl{For a finitely
generated group $G$, any irreducible deformation space which contains
a $G$--tree with finitely generated vertex groups is contractible.}
\medskip

The outline of the proof is as follows: starting with an unnormalized
deformation space \GT, we look at the space \GM(\GT) of morphisms
between elements of \GT.  A \textit{morphism} is a $G$--equivariant map
between $G$--trees which on each segment either folds or is an
isometry.  Given a morphism $\phi\co T \to Y$ we show that we can
continuously interpolate between the two $G$--trees.  We then fix some
reduced $G$--tree $T \in \GT$.  For another $G$--tree $Y \in \GT$ we
define a map $B(Y)\co T \to Y$, which is not a morphism but is nice in
certain respects.  The assignment $Y \mapsto B(Y)$ is a continuous
function between the appropriate spaces.  We redefine the metric
on $T$ to obtain another $G$--tree $T_{Y}$ (equivariantly homeomorphic
to $T$) such that $B(T)\co T_{Y} \to Y$ is a morphism.  Thus we can
homotope \GT \ to the space of trees equivariantly homeomorphic to
$T$.  We show this space is homeomorphic to a simplex$\times \R^+$,
thus $\GT$ is contractible.  This homotopy descends to a contraction of
the deformation space \GD.

Originally, the following proof was only for finitely generated
generalized Baum-slag--Solitar groups, for which there is a natural
$\Out(G)$--invariant deformation space.  A \textit{generalized
Baumslag--Solitar group} is a group which admits an action on a
simplicial tree where the stabilizer of any point is isomorphic to \Z.
However after a research announcement by Guirardel and Levitt
\cite{ar:GL}, which contains Theorem \ref{th:d_contractible}, we
noticed that our proof for generalized Baumslag--Solitar groups went
through in the general case after modifying case (ii) in Lemma
\ref{lm:step1}.  We are grateful for their announcement.  They have
proven Theorem \ref{th:d_contractible} in the case of a free product
and have given several consequences \cite{ar:GL2}.

The majority of material presented within is in Skora's preprint
\cite{ar:S}.  As this preprint was never published, we present the
full details here.  The main difference from \cite{ar:S} is section
\ref{sc:cont_def}.

\textbf{Acknowledgements}\qua  This work was done under the supervision of
my advisor Mladen Bestvina.  In addition to thanking him for the
helpful discussions, I am also grateful for discussions with Lars
Louder and for the research announcement of Vincent Guirardel and
Gilbert Levitt.  Thanks are also due to the referee for suggestions
improving the exposition.

%%%%%%%%%%%%%%%%%%%%%%%%%%%%%%%%%%%%%%%%%%%%%%%%%%%%%%%%%%%%%%%%%%%%%%%%%%%%% 

\section{Preliminaries}\label{sc:pre}

For a $G$--tree $T$, the length function $l_T\co G \to [0,\infty)$ is
defined by $l_T (g) = \min_{x \in T} d(x,gx)$.  The
\textit{characteristic set} $T_g$, of a element $g \in G$ is where
this minimum is realized, i.e.\ $T_g = \{ x \in T \ | \ d(x,gx) =
l_T(g) \}$.  An element is \textit{elliptic} if $l_T(g)=0$ and
\textit{hyperbolic} otherwise.  For $g \in G$ hyperbolic, the
characteristic set is isometric to \R \ and $g$ acts on $T_g$ by
translation by $l_T(g)$.  In this case the characteristic set of $g$
is often called the \textit{axis} of $g$. Note that $d(x,gx) =
2d(x,T_g) + l_T(g)$ for both $g$ elliptic or $g$ hyperbolic.  If a
subgroup $H \subseteq G$ is elliptic, we define the characteristic set
of $H$ as $T_H = \{ x \in T \ | \ hx = x \ \forall h \in H \}$.  For a
closed set $A \subseteq T$ we let $p_A\co T \to A$ denote the nearest
point projection.  A map between metric simplicial trees $\phi\co T
\to T'$ is morphism if for any segment $[x,y] \subseteq T$ there is a
subsegment $[x,x'] \subseteq [x,y]$ on which $\phi$ is an isometry.
If $T$ and $T'$ are $G$--trees, we also require that $\phi$ is
$G$--equivariant.

We have following dictionary of group actions on trees \cite{ar:CM}.
A $G$--tree is \textit{trivial} if there is a fixed point and
\textit{minimal} if there is no proper invariant subtree.  A $G$--tree $T$
is \textit{reducible} if:

\begin{itemize}

\item[(1)] every element fixes a point (equivalent to being trivial for 
finitely generated groups); or

\item[(2)] $G$ fixes exactly one end of $T$; or

\item[(3)] $G$ leaves a set of two ends of $T$ invariant.

\end{itemize}

\noindent
If $T$ is not reducible, it is \textit{irreducible}.  A $G$--tree is
irreducible if and only if there are two hyperbolic elements whose
axes are either disjoint or intersect in a compact set \cite{ar:CM}.
This feature is preserved by elementary deformations \cite{ar:F1},
hence any $G$--tree obtained via an elementary deformation from an
irreducible $G$--tree is also irreducible.

Unless otherwise stated, we will always assume $G$--trees are minimal
and irreducible.  We say a deformation space is \textit{irreducible}
if every $G$--tree in the space is irreducible.  By the above
statement, a deformation space is irreducible if any $G$--tree in the
space is irreducible.

%%%%%%%%%%%%%%%%%%%%%%%%%%%%%%%%%%%%%%%%%%%%%%%%%%%%%%%%%%%%%%%%%%%%%%%%%%%%% 

\section{Topology on deformation spaces}\label{sc:top}

We endow an unnormalized deformation space \GT \ with the
\textit{Hausdorff--Gromov} topology.  Gromov introduced this topology
as a way to compare two distinct metric spaces \cite{ar:G}.  This
topology generalizes the Hausdorff distance between two closed sets in
a metric space.  The deformation space \GD \ is then topologized as
the quotient $\GT / \R^+$.

The Hausdorff--Gromov topology is defined as follows.  Let $X,Y$ be
metric $G$--spaces, i.e.\ metric spaces equipped with isometric
$G$--actions.  For any $\epsilon>0$, an
$\epsilon$--\textit{approximation} is a set $R \subseteq X \times Y$
that surjects onto each factor such that if $x,x' \in X$ and $y,y' \in
Y$ with $xRy$ (i.e.\ $(x,y) \in R$) and $x'Ry'$ then $|d(x,x') -
d(y,y')|<\epsilon$.  We say that $R$ is a \textit{closed}
$\epsilon$--approximation if $R$ is closed in $X \times Y$.  For a
finite subset $P \subseteq G$ and subspaces $K \subseteq X,L \subseteq
Y$ the $\epsilon$--approximation in $K \times L$ is
$P$--\textit{equivariant} if whenever $g \in P$, $x,gx \in K$ and $y
\in L$ with $xRy$ then $gy \in L$ and $gxRgy$.

Given an $\epsilon$--approximation $R \subseteq X \times Y$, we let
$R_\delta$ denote the closed $\delta$--neighborhood of $R$ using the
$L^1$ metric.  In other words $R_\delta = \{(x,y) \in X \times Y \ | \
\exists (x',y') \in R \mbox{ with } d(x,x') + d(y,y') \leq
\delta\}$. One can show that $R_{\delta}$ is a $(\epsilon +
2\delta)$--approximation.  If $R$ is $P$--equivariant, then
$R_{\delta}$ is $P$--equivariant.

These $\epsilon$--approximations can topologize any set of metric
$G$--spaces.  In particular, they can topologize any unnormalized
deformation space \GT.  Let \GSs \ be such a set of metric
$G$--spaces.  Then for $X \in \GSs$, $K \subseteq X$ compact, $P
\subseteq G$ finite and $\epsilon>0$ define a basic open set
$U(X,K,P,\epsilon)$ to be the set of all $Y \in \GSs$ such that there
is a compact set $L \subseteq Y$ and a $P$--equivariant closed
$\epsilon$--approximation $R \subseteq K \times L$.  If $K \subseteq
K'$ and $P \subseteq P'$ then $U(X,K',P',\epsilon) \subseteq
U(X,K,P,\epsilon)$.  This will allow us to assume that certain subsets
of $X$ and $G$ are contained in $K$ and $P$ respectively by shrinking
our basic open set.

Given an $\epsilon$--approximation $R \subseteq X \times Y$, we will
assume it is \textit{full}: i.e.\ if $xRy$ and $x'Ry'$ then every point
in $[x,x']$ is related by $R$ to some point in $[y,y']$ and vice
versa.  This is not necessary but it cleans up some of the proofs in
sections \ref{sc:basepoint} and \ref{sc:cont_def}.  When the set
$\GSs$ contains only trees the two topologies generated are the same.
For $X,P,\epsilon$ as above let $U_f(X,K,P,\epsilon)$ be the set of
all $Y \in \GSs$ such that there is a finite subtree $L \subseteq Y$
and a $P$--equivariant closed full $\epsilon$--approximation $R
\subseteq K \times L$.  Clearly we have $U_f(X,K,P,\epsilon) \subseteq
U(X,K,P,\epsilon)$.  We now show the opposite inclusion of bases.

For two trees $X,Y,$ subsets $K \subseteq X,L \subseteq Y$ related by
an $\epsilon$--approximation $R \subseteq K \times L$ and a finite
segment $[x_1,x_2] \subseteq K$, let $R([x_1,x_2]) = \{z \in L \ | \ \exists
x \in [x_1,x_2] \mbox{ with } xRz \}$.  For $z \in R([x_1,x_2])$
with $x_i R y_i$ for some $y_i \in Y, i =1,2$ we have $d(z,[y_1,y_2])
< 2 \epsilon$.  We have the following statement about the density of
$R([x_1,x_2])$.
 
\begin{lemma}\label{lm:Za_dense}
If $z_0 \in [y_1,y_2] \subseteq Y$ where $x_i R y_i$ for $i=1,2$, then there 
is a $z \in R([x_1,x_2])$ such that $d(z_0,z) < 2 \epsilon$.
\end{lemma}

\begin{proof}
We assume this is not the case.  Let $d(y_1,z_0) = d_1, d(y_2,z_0) =
d_2$.  As $y_1,y_2 \in R([x_1,x_2])$ we can assume both $d_1$ and
$d_2$ are larger than $\epsilon$.  Take $x \in [x_1,x_2]$ such that
$d(x_1,x) = d_1$.  Therefore $d(x_2,x) < d_2 + \epsilon$.  There is a
$z \in R([x_1,x_2])$ such that $xRz$.  For this $z$, $d(y_1,z) < d_1 +
2 \epsilon$ and $d(y_2,z) < d_2 + 2 \epsilon$.  Now we let $z' =
p_{[y_1,y_2]}(z)$.  Hence by our initial assumption $d(z_0,z') +
d(z,z') \geq d(z_0,z) \geq 2 \epsilon$.  Assume without loss of
generality that $z'$ is closer to $y_1$ than $z_0$ is.  Then $d(y_2,z)
= d(y_2,z_0) + d(z_0,z') + d(z',z) \geq d_2 + 2 \epsilon$, a
contradiction.
\end{proof}

To finish up the claim that the two above mentioned topologies are the
same we let $\delta = \frac{\epsilon}{5}$.  Then for $Y \in
U(X,K,P,\delta)$ we have a $P$--equivariant $\delta$--approximation
between $K$ and some finite subtree $L \subseteq Y$. By the above
$R_{2\delta}$ is a full $P$--equivariant $\epsilon$--approximation
between $K$ and $L$.  Therefore $U(X,K,P,\delta) \subseteq
U_f(X,K,P,\epsilon)$ and the two topologies are indeed the same.

We will also topologize the space of morphisms between elements in a
deformation space.  Let $\phi\co X \to X'$ and $\psi\co Y \to Y'$ be
$G$--equivariant maps, a closed $\epsilon$--approximation between these two
maps is a pair $(R,R')$ such that: 

\begin{itemize}

\item[(1)] $R \subseteq X \times Y$ and $R' \subseteq X' \times Y'$ are
closed $\epsilon$--approximations; and

\item[(2)] for $x \in X,y \in Y$ if $xRy$ then $\phi(x)R'\psi(y)$.

\end{itemize}

\noindent
Let $P \subseteq G$ be finite and $K \subseteq X, K' \subseteq X', L
\subseteq Y, L' \subseteq Y'$ be subspaces.  The
$\epsilon$--approximation $(R,R')$ is $P$--equivariant if $R$ and $R'$
are $P$--equivariant on the appropriate subspaces.  Note that if $z \in
\gr(\phi)$ then by definition there is a $w \in \gr(\psi)$ with
$z(R,R')w$.

As above this allows us to topologize a set of $G$--equivariant maps
between $G$--spaces.  In particular we can topologize \GM(\GT), the
set of morphisms between elements of \GT.  Let $\GMs$ a set of
$G$--equivariant maps between $G$--spaces.  For $\phi\co X \to X'$ in
$\GMs$, and $K \subseteq X, K' \subseteq X'$ both compact with
$\phi(K) \subseteq K'$, $P \subseteq G$ finite and $\epsilon>0$ define
the basic open set $U(\phi,K \times K',P,\epsilon)$ to be the set of
all maps $\psi\co Y \to Y'$ in $\GMs$ such that there are compact sets
$L \subseteq Y$, $L' \subseteq Y'$ with $\psi(L) \subseteq L'$ and a
$P$--equivariant closed $\epsilon$--approximation $(R,R')$ bewtween
$\phi\co K \to K'$ and $\psi\co L \to L'$.

For a space $\GSs$ of metric $G$--spaces and a space $\GMs$ of
$G$--equivariant maps between the elements of $\GSs$ we have the two
continuous maps $\Do$ and $\Ra$ defined from $\GMs$ to $\GSs$ which send a
map to its domain and range respectively.  In other words, for $\phi\co X
\to X'$ an element of $S'$ we have $\Do (\phi) = X$ and $\Ra (\phi) = X'$.

There are two other topologies one might use to topologize a
deformation space.  Let \CC \ be the set of conjugacy classes for
$G$.  Then we have a function $l\co \GT \to \R^{\CC}$ where the
coordinates are given by the length functions $l_T(c)$ where $c \in
\CC$.  Culler and Morgan showed that for minimal irreducible actions
on \R-trees this function is injective \cite{ar:CM}.  This defines a
topology on \GT \ (and hence on \GD) called the \textit{axes}
topology.  Paulin proved that for spaces of minimal irreducible
actions on \R-trees, the Hausdorff--Gromov topology is the same as the
axes topology \cite{ar:P}.

We can define the \textit{weak} topology directly on \GD.  The
\textit{volume} of a $G$--tree $T$, denoted $\vol(T)$, is the sum of
the lengths of the unoriented edges of $T/G$.  We identify \GD \
with the $G$--trees in \GT \ that have volume one.  By reassigning the
lengths of the edges of $T/G$ in a manner to hold the volume constant
we can define a simplex in \GD.  The weak topology is defined by
considering \GD \ as the union of such simplicies.  In general, the
weak topology is different from the axes and Hausdorff--Gromov
topology, see \cite{ar:MM} for an example.

%%%%%%%%%%%%%%%%%%%%%%%%%%%%%%%%%%%%%%%%%%%%%%%%%%%%%%%%%%%%%%%%%%%%%%%%%%%%% 

\section{Deforming trees}\label{sc:def_tree}

A morphism $\phi\co T \to T'$ between trees in an unnormalized
deformation space \GT \ can be decomposed into elementary deformations
\cite{ar:F1}.  We will define trees $T_t$ which continuously
interpolate between $T$ and $T'$.

For the morphism $\phi\co T \to T'$, a nontrivial segment
$[x,x'] \subseteq T$ is \textit{folded} if $\phi(x)=\phi(x')$.  A
folded segment is \textit{maximally folded} if it cannot be locally
extended to a segment which is folded.  On a maximally folded segment
$[x,x']$ the function $d(\phi(z),\phi(x))$ attains a local maximum at
possibly several points.  Such points are called \textit{fold points}
of the morphism $\phi$.  The points at where the global maxima are
obtained are called \textit{maximal fold points}.  We remark that
every fold point is a maximal fold point for some maximally folded
segment.  A fold point $z$ is $d$-\textit{deep} if
$d(\phi(x),\phi(z))>d$ for some maximally folded segment $[x,x']$ of
which $z$ is a maximal fold point.

\begin{figure}[ht!]\small\anchor{fg:thicken}
\setlength{\unitlength}{0.0003in}
\gdef\SetFigFont#1#2#3#4#5{\relax}%
\begin{center}
{\renewcommand{\dashlinestretch}{30}
\begin{picture}(13104,3841)(0,-10)
\path(1680,3102)(987,3795)  %left top v-l
\path(1680,3102)(2382,3795) %left top v-r
\path(12,3795)(987,3795)  %left top left
\path(2382,3795)(3357,3795) %left top right
\path(1680,987)(987,1654) %left bottom y-l
\path(1680,987)(2382,1654) %left bottom y-r
\path(1680,987)(1680,12)  %left bottom y-b
\path(6627,2565)(5832,1530)  %center left y-r
\path(5043,2576)(5838,1541)(5838,176) %center left y-l&b
\path(5843,176)(8622,1350) %center bottom
\path(6627,2565)(9417,3750) %center right top 
\path(5043,2576)(7829,3755) %center left top
\path(7829,3755)(8624,2720) %center right y-l
\path(9417,3750)(8622,2715)(8622,1350) %center right y-r&b
\path(5835,1537)(8610,2707) %center midline
\path(5043,2576)(6672,1905)(7287,780)(7512,2220)(9417,3750) %center graph
\path(5337,2205)(6417,1755)(7025,685) %center left inside
\path(5367,2715)(7002,1995)(7152,1740)
	(7242,2115)(9090,3605) %center center inside
\path(7628,928)(7842,2355)(9110,3340) %center right inside
\path(10707,2820)(11682,2820) %right left
\path(11682,2820)(11907,2595) %right y-l
\path(12117,2820)(11907,2595)(11907,2070) %right y-r&b
\path(12117,2820)(13092,2820) %right right
\path(1677,2820)(1677,1890) %left arrow
\blacken\thicklines
\path(1639.500,2010.000)(1677.000,1890.000)(1714.500,2010.000)(1677.000,1974.000)(1639.500,2010.000)
\thinlines
\put(1212,2310){\makebox(0,0)[lb]{\smash{{{\SetFigFont{12}{24.0}{\rmdefault}{\mddefault}{\updefault}$\phi$}}}}}
\put(7662,510){\makebox(0,0)[lb]{\smash{{{\SetFigFont{12}{24.0}{\familydefault}{\mddefault}{\updefault}$W_t$}}}}}
\put(11772,1410){\makebox(0,0)[lb]{\smash{{{\SetFigFont{12}{24.0}{\familydefault}{\mddefault}{\updefault}$T_t$}}}}}
\end{picture}}
\end{center}
\caption{$W_t$ and $T_t$ for the morphism on the left}\label{fg:thicken}
\end{figure}
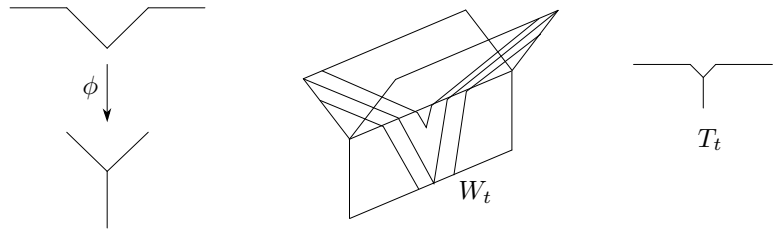

We let $m(\phi) = \sup \{ d(\phi(z),\phi(x)) \ | \ z \in [x,x'] \mbox{
where } \phi(x)=\phi(x') \}$.  Then $m(\phi)$ is finite as elementary
deformations are quasi-isometries \cite{ar:F1}.  Notice that $m(\phi)
= 0$ if and only if $\phi$ is an isometry and hence $T=T'$ as
$G$--trees.  For $0 \leq t \leq 1$ we define $V_t = \{ (x,y) \in T
\times T' \ | \ d(\phi(x),y) \leq m(\phi)t \}$.  For $(x,y) \in V_t$
let $C_t(x,y)$ denote the path component of $V_t \cap (T \times
\{y\})$ which contains $(x,y)$.  Finally, we define:
$$W_t = \{ (x,y) \in V_t \ | \ C_t(x,y) \cap \gr(\phi) \neq \emptyset
 \}.$$
Thus $W_t$ is a thickening of $\gr(\phi) \subseteq T \times T'$.  We
will write $W_t(\phi)$ when we need to specify the morphism.  Let
$\F_t$ be a partition of $W_t$ into sets which are the path components
of $W_t \cap (T \times \{ y \})$ for $y \in T$ and $T_t = W_t /
\F_t$.  We denote points in $T_t$ by $[z]_t$ for $z \in W_t$.  As
$\F_t$ is $G$--equivariant, $T_t$ is a $G$--tree.  For $0 \leq s \leq t
\leq 1$ the inclusions $W_s \to W_t$ induce $G$--equivariant maps
$\phi_{st}\co T_s \to T_t$, \figref{fg:thicken}.

\begin{figure}[ht!]\small\anchor{fg:tcp}
\setlength{\unitlength}{0.000275in}
\begin{center}
{\renewcommand{\dashlinestretch}{30}
\begin{picture}(7734,7050)(0,-10)
\texture{55888888 88555555 5522a222 a2555555 55888888 88555555 552a2a2a 2a555555 
	55888888 88555555 55a222a2 22555555 55888888 88555555 552a2a2a 2a555555 
	55888888 88555555 5522a222 a2555555 55888888 88555555 552a2a2a 2a555555 
	55888888 88555555 55a222a2 22555555 55888888 88555555 552a2a2a 2a555555 }
\put(6357,3318){\shade\ellipse{130}{130}}
\put(6357,3318){\ellipse{130}{130}}
\put(1113,3882){\shade\ellipse{130}{130}}
\put(1113,3882){\ellipse{130}{130}}
\path(12,5448)(1572,33)(2907,33)(3867,4353)
\path(2237,30)(632,5580)
\path(2247,18)(3807,7008)(3807,7023)
\path(1317,5583)(2157,2523)(3147,7023)(4347,7023)
\path(5402,13)(3842,7003)(3842,7018)
\path(5416,28)(7021,5578)
\path(7722,5427)(6162,12)(4827,12)(3867,4332)
\path(6287,5580)(5447,2520)(4457,7020)(3257,7020)
\thicklines
\path(1137,3828)(1542,2508)(2802,2508)
	(3237,4353)(4452,4353)(4872,2478)
	(6117,2478)(6342,3273)
\put(1137,4158){\makebox(0,0)[lb]{\smash{{{\SetFigFont{12}{14.4}{\rmdefault}{\mddefault}{\updefault}$z_1$}}}}}
\put(6627,3483){\makebox(0,0)[lb]{\smash{{{\SetFigFont{12}{14.4}{\rmdefault}{\mddefault}{\updefault}$z_2$}}}}}
\end{picture}}
\end{center}
\caption{A taut corner path in the subset $W_t$ between the points $z_1$ and $z_2$.  The central line is the graph of the morphism.  The essential pieces are the segments which lie in the graph; the nonessential pieces are the horizontal segments.}\label{fg:tcp}
\end{figure}
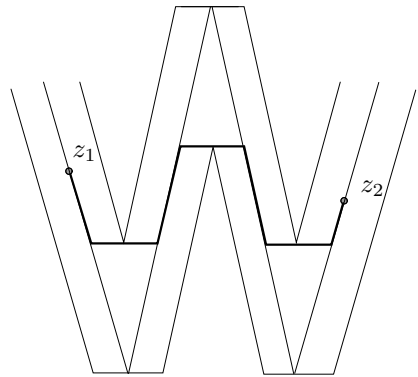

A path $\gamma:[0,1] \to W_t$ is \textit{taut} if for components $A$
in $\F_t$, $\gamma^{-1}(A)$ is connected.  For $t>0$, a
non-backtracking path $\gamma$ in $T_t$ lifts to a path
$\tilde{\gamma}$ in $W_t$ with endpoints in $\gr(\phi)$. This lift
$\tilde{\gamma}$ is homotopic relative to these endpoints to a taut
product of paths $\gamma_1 \cdots \gamma_k$ where each $\gamma_i$ lies
either in a component of $\F_t$ or is a non-backtracking path in
$\gr(\phi)$.  For $t=0$ a non-backtracking path $\gamma$ in $T_0$
lifts to a path $\tilde{\gamma}$ which is homotopic relative to its
endpoints to a taut product of paths $\gamma_1 \cdots \gamma_k$ where
$\phi$ is an isometry on each $\gamma_i$.  We call these
decompositions \textit{taut corner paths}, the pieces lying in
$\gr(\phi)$ are called \textit{essential}, the pieces lying in some
component of $\F_t$ are called \textit{nonessential}, see
\figref{fg:tcp}.  Metrize $T_t$ by setting $\length(\gamma)$ equal to
the sum of the lengths of the essential pieces measured in $T'$ (or
equivalently measured in $T$).  With this metric the maps $\phi_{st}$
are morphisms.

\begin{lemma}\label{lm:def_phi}
For the above definitions: $T_0 = T$, $T_{1} = T'$ as $G$--trees,
$\phi_{00} = Id_{T}$ and $\phi_{01} = \phi$.
\end{lemma}

\begin{proof}
The only nonobvious claim here is $T_1 = T'$.  This is equivalent to
saying that the sets $W_1 \cap (T \times \{y \})$ are connected.  Let
$(x_1,y),(x_2,y) \in W_1 \cap (T \times \{y \})$.  We will show that
these two points lie in the same component.  Choose $(z_1,y), (z_2,y)
\in W_1 \cap (T \times \{ y \})$ such that $\phi(z_i) = y$ and
$(x_i,y),(z_i,y)$ are in the same component of $W_1 \cap (T \times \{
y\})$ for $i = 1,2$.  For $z \in [z_1,z_2]$ we have $d(\phi(z),y) \leq
m(\phi)$.  Thus the pairs of points $(z_1,y),(z_2,y)$ are in the same
component of $W_1 \cap (T \times \{ y \})$.  Then as $(x_1,y)$ is in
the same component as $(z_1,y)$ and $(x_2,y)$ is in the same component
as $(z_2,y)$, the points $(x_1,y)$ and $(x_2,y)$ are in the same
component.  Thus $W_1 \cap (T \times \{y \})$ is connected.
\end{proof}

\begin{lemma}\label{lm:irreducible}
If $T'$ is irreducible, then so is $T_t$ for $0 \leq t \leq 1$.
\end{lemma}

\begin{proof}
As the $G$--tree $T'$ is irreducible, there are $g,h \in G$ which act
hyperbolically on $T'$ such that $T_g \cap T_h$ is empty or compact
\cite{ar:CM}.  As equivariant maps cannot make elliptic elements act
hyperbolically, $g,h$ act hyperbolically in $T_t$.  The maps
$\phi_{st}$ are quasi-isometries, hence the axes of $g$ and $h$ have
empty or compact intersection.  This implies that the $G$--tree $T_t$
is irreducible.
\end{proof}

The following lemma is obvious.

\begin{lemma}\label{lm:def_space}
$T_t$ is in the same unnormalized deformation space as $T$ and $T'$
for $0 \leq t \leq 1$.
\end{lemma}

\begin{remark}\label{rm:R-invar}{\rm
For future reference we remark that the above construction is
invariant under the $\R^+$--action.  In other words if we scale both
$T$ and $T'$ by a nonzero positive number $k$, then the trees $T_t$
are scaled by $k$.}
\end{remark}

%%%%%%%%%%%%%%%%%%%%%%%%%%%%%%%%%%%%%%%%%%%%%%%%%%%%%%%%%%%%%%%%%%%%%%%%%%%%% 

\section{Continuity of deformation}\label{sc:con_def}

Fix an unnormalized deformation space \GT.  Recall that \GM(\GT) is
the space of all morphisms between $G$--trees in \GT.  Define $\Phi\co
\GM (\GT) \times \{ (s,t) \ | \ 0 \leq s \leq t \leq 1 \} \to \GM (\GT)$
by $\Phi(\phi,(s,t)) = \phi_{st}$.  The goal of this section is the
following theorem:

\begin{theorem}\label{th:Phi}
$\Phi$ is continuous.
\end{theorem}

We have some work before we can prove this.  The approach is the same
as in Skora's preprint \cite{ar:S}, with the addition of Lemmas
\ref{lm:close_m} and \ref{lm:nonessential}.  We will consider a fixed
morphism $\phi\co X \to X'$ between finite simplicial trees and prove
some results about morphisms $\psi\co Y \to Y'$ which are close to
$\phi$.  The main step is to show Lemma \ref{lm:close_length}: if the
two morphisms $\phi\co X \to X'$ and $\psi\co Y \to Y'$ are close and
we fold both $X$ and $Y$ for a similar amount of time, then the two
folded trees have comparable lengths.  To prove this, we show that for
a taut corner path in $W_t(\phi)$, the individual pieces are related
to taut corner paths of comparable length in $W_s(\psi)$ when both $s$
and $t$ are close and $\phi$ and $\psi$ are close.

Our first step is to show that maps close to $\phi$ have similar
folding data.

\begin{lemma}\label{lm:close_m}
Let $\phi\co X \to X'$ be a morphism of finite simplicial trees.  Then
for all $\epsilon>0$ there is a $\delta>0$ such if $(R,R')$ is a
$\delta$--approximation between $\phi\co X \to X'$ and $\psi\co Y \to Y'$
then $|m(\phi) - m(\psi)|<\epsilon$.
\end{lemma}

\begin{proof}
Let $\epsilon>0$ be arbitrary.  There are two cases to deal with.
Notice that the lemma is symmetric with respect to $\phi$ and $\psi$.

\textbf{Case 1}\qua $m(\phi),m(\psi)>0$

Set $\delta = \min\{ \frac{\epsilon}{2},\frac{m(\phi)}{3},
\frac{m(\psi)}{3} \}$.  Let $z \in [x,x']$ be such that
$d(\phi(z),\phi(x)) = m(\phi)$ and $\phi(x)=\phi(x')$.  There are
corresponding points $y,y',w \in Y$ such that $xRy,x'Ry'$ and $zRw$.
As we can assume that $R$ is full, we may assume that $w \in [y,y']$.
Then $d(\psi(y),\psi(y'))< \delta$ and
$d(\psi(w),\psi(y)),d(\psi(w),\psi(y')) > m(\phi) - \delta > \delta$.
Hence there is a subsegment contained in $[y,y']$ and containing $w$
which is folded.  Thus $m(\psi) > m(\phi) - 2\delta$.  Repeating
the argument for $Y$ we see that $m(\phi) > m(\psi) - 2\delta$.
Hence we see that $|m(\phi) - m(\psi)| < \epsilon$.

\textbf{Case 2}\qua $m(\phi)=0$ and $m(\psi)>0$ 

Set $\delta = \frac{\epsilon}{2}$.  Let $[y,y'] \subseteq Y$ be a
folded segment where $w \in [y,y']$ attains $d(\psi(w),\psi(y)) =
m(\psi)$.  For corresponding points $x,x',z \in X$, we have that
$[x,x']$ is embedded and $d(\phi(x),\phi(x')) < \delta$.  Hence for
all $z \in [x,x']$, we have $d(\phi(z), \phi(x)) < \delta$.  Thus
$m(\psi) < \epsilon$.
\end{proof}

Let $F(\phi)$ denote the number of fold points for the morphism $\phi\co
X \to X'$.  Thus for $N(\phi) = 3(F(\phi)+1)$ we have that any
taut corner path $\gamma$ in $W_t(\phi)$ can be written as $\gamma =
\gamma_1 \cdots \gamma_n$ with $n \leq N(\phi)$ where each $\gamma_i$
is either essential or nonessential.

We need a similar statement about morphisms close to $\phi$.  It is
easy to see that we cannot expect a universal bound, but we can bound
the number of large folds, which is sufficient.  For $d > 0 $, we
introduce an equivalence relation on the set of fold points defined by
$z \sim_d z'$ if there is a sequence of fold points: $z = z_0,
\ldots,z_n = z'$ such that $d(z_i,z_{i+1}) < 2d$.  Let

\begin{center}
$F_d(\psi) = \left\{ \{ \mbox{fold points for }\psi \} / \sim_d
\right\} \setminus \{ \mbox{classes without a $d$-deep
point} \}$.
\end{center}

Notice that $|F_{m(\psi)s}(\psi)|$ is the number of fold points for
the map $\psi_{s1}\co Y_s \to Y'$.  Suppose $(R,R')$ is a
$\delta$--approximation between $\phi\co X \to X'$ and $\psi\co Y \to
Y'$ where $\delta \leq d$.  Then $|F_d(\psi)|$ is bounded independent
of $\psi$ as for each class in $F_d(\psi)$ we have a $\frac{d}{2}$
neighborhood in $X$, and the neighborhoods for different classes are
disjoint.  Set $F_d(\phi)$ to be the maximum of $|F_d(\psi)|$ over all
morphisms $\psi\co Y \to Y'$ for which there is a $d$--approximation
between $\phi\co X \to X'$ and $\psi\co Y \to Y'$.  As above we let
$N_d(\phi) = 3(F_d(\phi) + 1)$.  Thus if $\zeta$ is a taut corner path
in $W_s(\psi)$ then we can write $\zeta = \zeta_1 \cdots \zeta_n$ with
$n \leq N_d(\phi)$ where each $\zeta_i$ is either nonessential or has
length equal to the length of its image in $Y'$.

We now show that for a taut corner path in $W_t(\phi)$, the individual
pieces are related to a taut corner path in $W_s(\psi)$ of comparable
length.  This is proven for the essential pieces first.  As a
convention when taking several points in $X$ and points related to
them in $Y$, if some of the points in $X$ are the same we require
that the related points in $Y$ are the same.

\begin{lemma}\label{lm:essential}
Let $\phi\co X \to X'$ be a morphism of finite simplicial trees.  Let
$z_1,z_2$ be points in $\gr(\phi) \subseteq W_t(\phi)$ such that the
taut corner path $\gamma$ between them lies entirely in $\gr(\phi)$.
Then for all $\epsilon>0$ there is a $\delta>0$ such that if $(R,R')$
is a $\delta$--approximation between $\phi\co X \to X'$ and $\psi\co Y
\to Y'$ and $w_i \in \gr(\psi)$ where $z_i(R,R')w_i$ for $i=1,2$, then
$|\length(\gamma)-\length(\zeta)| < \epsilon$ where $\zeta$ is the taut
corner path $\zeta$ in $W_s(\psi)$ from $w_1$ to $w_2$,
\end{lemma}

\begin{proof}
Let $\epsilon>0$ be arbitrary.  Let $\delta = \epsilon$ and assume the
data in the hypothesis.  Let $z_i=(x_i,x_i'),w_i=(y_i,y_i')$ for
$i=1,2$ and let $\zeta$ be the taut corner path in $W_s(\psi)$
connecting $w_1$ to $w_2$.  By hypothesis $\length(\gamma) = d(x_1,x_2)
= d(x_1',x_2')$.  As $\length(\zeta) \leq d(y_1,y_2) < d(x_1,x_2) +
\delta$ and $\length(\zeta) \geq d(y_1',y_2') > d(x_1',x_2') - \delta$,
we have the conclusion of the lemma.
\end{proof}

Next we have a similar statement for the nonessential pieces:

\begin{lemma}\label{lm:nonessential}
Let $\phi\co X \to X'$ be a morphism of finite simplicial trees.  Let
$z_1,z_2$ be points in $\gr(\phi) \subseteq W_t(\phi)$ such that the
taut corner path $\gamma$ between them lies entirely in a component of
$\F_t$.  Then for all $\epsilon>0$ there is a $\delta>0$ such that if
$(R,R')$ is a $\delta$--approximation between $\phi\co X \to X'$ and
$\psi\co Y \to Y'$, $|m(\psi)s-m(\phi)t|< \delta$, and $w_i \in
\gr(\psi)$ where $z_i(R,R')w_i$ $i=1,2$, then $\length(\zeta) <
\epsilon$ where $\zeta$ is the taut corner path in $W_s(\psi)$ from
$w_1$ to $w_2$,
\end{lemma}

\begin{proof}
Let $\epsilon>0$ be arbitrary.  We have two cases depending on
$m(\phi)$ and $t$.  Let $(R,R')$ be a $\delta$--approximation with
$z_i,w_i$ as in the statement above where $\delta$ is chosen in the
individual cases.  Say $z_i = (x_i,x_i'),w_i = (y_i,y_i')$ for
$i=1,2$.  Then from the definitions we have $x_1' = x_2'$ and
$\phi([x_1,x_2])$ stays within $m(\phi)t$ of $x_1'$.

\textbf{Case 1}\qua $t=0$ or $m(\phi) = 0$

Let $\delta = 1$.  Then as $z_1=z_2$, we have $w_1=w_2$ by the
above convention.  Hence $\length(\zeta)=0$.

\textbf{Case 2}\qua $t>0$ and $m(\phi) > 0$

Let $N = N_d(\phi)$ as above where $d = \frac{2m(\phi)t}{3}$ and set
$\delta = \min\{\frac{\epsilon}{2N},\frac{m(\phi)t}{3}\}$.  As
$m(\psi)s > m(\phi)t - \delta \geq d$, the number of fold points for
$\psi_{s1}$ is less than $F_d(\phi)$.  Therefore we can write $\zeta =
\zeta_1 \cdots \zeta_n$ where $n \leq N$ and each $\zeta_i$ is
nonessential or has length equal to the length of its image in $Y'$.

If $\psi([y_1,y_2])$ is contained within a $m(\psi)s$ neighborhood
about $y_1'$, then $\zeta$ is nonessential.  This might not be the
case, but the length of an essential piece of $\zeta$ is bounded by
how far $\psi([y_1,y_2])$ travels away from $y_1'$: $\length(\zeta_i)
\leq \max \{ \{ d(\psi(y),y_1') - m(\psi)s \ | \ y \in
[y_1,y_2]\},0\}.$ Now we use fullness of the approximations to see:
$\length(\zeta_i) \leq \max \{ \{ d(\phi(x),x_1') - m(\psi)s + \delta \
| \ x \in [x_1,x_2] \}, 0 \} \leq \max \{ \{ d(\phi(x),x_1') -
m(\phi)t + 2\delta \ | \ x \in [x_1,x_2] \}, 0 \} \leq 2\delta$.

Thus we have $\length(\zeta) \leq \sum \length(\zeta_i) \leq 2\delta n < 
\epsilon$. 
\end{proof}

Putting together the previous two lemmas we have:

\begin{lemma}\label{lm:tcp}
Let $\phi\co X \to X'$ be a morphism of finite simplicial trees.  Let
$z_1,z_2$ be points in $\gr(\phi \subseteq W_t(\phi)$ and $\gamma =
\gamma_1 \cdots \gamma_n$ the taut corner path between them.  Then for
all $\epsilon>0$ there is a $\delta>0$ such if $(R,R')$ is a
$\delta$--approximation between $\phi\co X \to X'$ and $\psi\co Y \to
Y'$, $|m(\psi)s-m(\phi)t|<\delta$ and $w_i \in \gr(\psi)$ where
$z_i(R,R')w_i$ for $i=1,2$, then there is a path $\zeta = \zeta_1
\cdots \zeta_n$ in $W_s(\psi)$ from $w_1$ to $w_2$ with each $\zeta_i$
a taut corner path which satisfies
$|\length(\gamma_i)-\length(\zeta_i)|<\epsilon$ for $i=1,\ldots,n$.
\end{lemma}

The next lemma is a converse to Lemma \ref{lm:nonessential} and the
proof is simpler as we know how many fold points $\phi$ has.  Recall
that the image of $z \in W_t$ in the quotient tree is denoted $[z]_t$.

\begin{lemma}\label{lm:con_nonessential}
Let $\phi\co X \to X'$ be a morphism of finite simplicial trees.  Then
for all $\epsilon >0$ there is a $\delta>0$ such that if $(R,R')$ is a
$\delta$--approximation between $\phi\co X \to X'$ and $\psi: Y \to
Y'$ and $|m(\psi)s-m(\phi)t|<\delta$ then $d([z_1]_t,[z_2]_t) <
\epsilon$ where $z_i \in \gr(\phi)$, $w_i \in \gr(\psi)$ with
$z_i(R,R')w_i$ for $i=1,2$ and $[w_1]_s = [w_2]_s$.
\end{lemma}

\begin{proof}
Let $\epsilon>0$ be arbitrary and $\delta = \frac{\epsilon}{2N}$,
where $N = N(\phi)$.  Let $\gamma = \gamma_1 \cdots \gamma_n$ be the
taut corner path from $z_1$ to $z_2$ where each piece is either
essential or nonessential and $n \leq N$.  Using the same argument as
in case 2 for \ref{lm:nonessential}, we can bound the lengths of the
$\gamma_i$ by $2\delta$.  Thus $\length(\gamma) \leq 2\delta n <
\epsilon$.
\end{proof}

Now using the previous two lemmas, we are able to show that close
morphisms which are folded for a similar amount of time have
comparable lengths.  We will also remove the dependence on the folding
data using Lemma \ref{lm:close_m}.

\begin{lemma}\label{lm:close_length}
Let $\phi\co X \to X'$ be a morphism of finite simplicial trees.  Then
for all $\epsilon >0$ there is a $\delta>0$ such that if $(R,R')$ is a
$\delta$--approximation between $\phi\co X \to X'$ and $\psi\co Y \to
Y'$ and $|s-t|<\delta$ then $|d([z_1]_t,[z_2]_t) - d([w_1]_s,[w_2]_s)|
< \epsilon$ where $z_i \in \gr(\phi)$, $w_i \in \gr(\psi)$ with
$z_i(R,R')w_i$ for $i=1,2$.
\end{lemma}

\begin{proof}
Let $\epsilon$ be arbitrary.  Set $\epsilon_1 = \epsilon_2 =
\frac{\epsilon}{4N}$, where $N=N(\phi)$.  Use these to find
$\delta_1,\delta_2$ from Lemma \ref{lm:tcp} and Lemma
\ref{lm:con_nonessential} respectively.  Let $\epsilon_3 = \frac{1}{2}
\min \{ \delta_1,\delta_2 \}$ and take $\delta_3$ from Lemma
\ref{lm:close_m} using $\epsilon_3$.  Finally set $\delta = \min\{
\epsilon_2,\frac{\epsilon_3}{m(\phi)},\delta_3 \}$.

The choice of these parameters implies that if $(R,R')$ is a
$\delta$--approximation between $\phi\co X \to X'$ and $\psi\co Y \to
Y'$, and $|s-t|<\delta$, then $|m(\psi) - m(\phi)| < \epsilon_3$.
Thus $m(\psi)s - m(\phi)t < (m(\phi) + \epsilon_3)s - m(\phi)t <
m(\phi)(s-t) + \epsilon_3 < \delta_1,\delta_2$ and similarly
$m(\phi)t-m(\psi)s < \delta_1,\delta_2$.  Therefore we can use Lemma
\ref{lm:tcp} and Lemma \ref{lm:con_nonessential}.

We can write the taut corner path connecting $z_1$ and $z_2$ as
$\gamma = \gamma_1 \cdots \gamma_n$ where $n \leq N$ and each
$\gamma_i$ is either essential or nonessential.  Hence by Lemma
\ref{lm:tcp} we have a path $\zeta = \zeta_1 \cdots \zeta_n$
connecting $w_1$ to $w_2$ where each piece is a taut corner path and
$|\length(\gamma_i) - \length(\zeta_i)|<\epsilon_1$.  Hence
$d([w_1]_s,[w_2]_s) \leq \sum \length(\zeta_i) < \sum \left(
\length(\gamma_i) + \epsilon_1 \right) < d([z_1]_t,[z_2]_t) +
\epsilon$.

If $d([w_1]_s,[w_2]_s) \leq d([z_1]_t,[z_2]_t) - \epsilon$, then as
$d([z_1]_t,[z_2]_t) = \sum \length(\gamma_i) < \left( \sum
\length(\zeta_i) \right) + \frac{\epsilon}{2}$ we get that
$d([w_1]_s,[w_2]_s) < \sum \length(\zeta_i) - \frac{\epsilon}{2}$.
Since the only folds of $[\zeta]$ in $Y_s$ are at the intersection
points of $[\zeta_i]$ with $[\zeta_{i+1}]$, there are two points $q_1$
and $q_2$ on $\zeta$ such that the length along $\zeta$ between these
two points is greater than $\frac{\epsilon}{2N}$ but these are the
same point in $Y_s$.  Thus for points $p_1,p_2 \in W_t(\phi)$ with
$p_i (R,R') q_i$ for $i=1,2$ we have $d([p_1]_t,[p_2]_t) >
\frac{\epsilon}{2N} - \delta \geq \epsilon_2$.  However the choice of
$\delta_2$ implies that $d([p_1]_t,[p_2]_t) < \epsilon_2$ by Lemma
\ref{lm:con_nonessential}.  Hence we have a contradiction.  Therefore
$|d([z_1]_t,z_2]_t) - d([w_1]_s,[w_2]_s)| < \epsilon$.
\end{proof}

Thus the folded trees have comparable lengths.  We can use this to
build an $\epsilon$--approximation between these trees.  For morphisms
$\phi\co X \to X'$, $\psi\co Y \to Y'$ which are related by an
$\epsilon$--approximation $(R,R')$ we define a new relation
$[R,R]_{ts}$ from $X_t$ to $Y_s$ by $[z]_t[R,R']_{ts}[w]_s$ whenever
$z(R,R')w$ for $z \in \gr(\phi)$ and $w \in \gr(\psi)$.  We now prove
a lemma about this relation when $s$ and $t$ are close.

\begin{lemma}\label{lm:new_relation}
Let $\phi\co X \to X'$ be a morphism of finite simplicial trees.  For
all $\epsilon>0$ there is a $\delta>0$ such that if $(R,R')$ is a
$\delta$--approximation between $\phi\co X \to X'$ and $\psi\co Y \to
Y'$ and $|s-t|<\delta$ then $[R,R']_{ts}$ is an
$\epsilon$--approximation from $X_t$ to $Y_s$.  If $R$ and $R'$ are
$P$--equivariant, then so is $[R,R']_{ts}$.
\end{lemma}

\begin{proof}
Let $\epsilon>0$ be arbitrary and choose $\delta$ from Lemma
\ref{lm:close_length}.  Given data as in the hypothesis, $[R,R']_{ts}$
is an $\epsilon$--approximation.  It also follows that if $R$ and $R'$
are $P$--equivariant, then so is $[R,R]_{ts}$.
\end{proof}

Given an arbitrary morphism $\phi\co T \to T'$ between $G$--trees in
the unnormalized deformation space \GT, for subtrees $X \subseteq T,X'
\subseteq T'$ such that $\phi(X) \subseteq X'$ we can define $X_t$ as
$W_t(\phi |_X) / (\F_t \cap (X \times X'))$.  We can now prove that
$\Phi$ is continuous.

\begin{proof}
Let $\epsilon>0$ be arbitrary.  Let $\phi\co T \to T'$ and $0 \leq s
\leq t \leq 1$ be given.  Assume $U$ is the basic open set around
$\phi_{st}$ given by $U=U(\phi_{st},X \times X',P,\epsilon)$ where $X
\subseteq T,X' \subseteq T'$ are finite subtrees and $P$ is a finite
subset of $G$.  Let $\delta$ be given by Lemma \ref{lm:new_relation},
and $V = U(\phi,X \times X',P,\delta)$.

Suppose $\psi\co\widetilde{T} \to \widetilde{T}'$ with $\psi \in V$ and
$|p-s|<\delta,|q-t|<\delta$.  We will show that $\psi_{pq} \in U$.

For some finite subtrees $Y \subseteq \widetilde{T}, Y' \subseteq
\widetilde{T}'$ there is a $\delta$--approximation $(R,R')$ from
$\phi:X \to X'$ to $\psi:Y \to Y'$.  The claim is that
$([R,R']_{sp},[R,R]_{tq})$ is a closed $\epsilon$--approximation from
$\phi_{st} : X_s \to X_t$ to $\psi_{pq}:Y_p \to Y_q$.  The choice of
$\delta$ implies that both $[R,R']_{sp}$ and $[R,R']_{qt}$ are
$\epsilon$--approximations by Lemma \ref{lm:new_relation}.  If $[z]_s
[R,R']_{sp} [w]_p$ then we have that $[z]_t [R,R']_{tq} [w]_q$.
Therefore $\psi_{pq} \in U$.
\end{proof}

%%%%%%%%%%%%%%%%%%%%%%%%%%%%%%%%%%%%%%%%%%%%%%%%%%%%%%%%%%%%%%%%%%%%%%%%%%%%% 

\section{Continuity of base point}\label{sc:basepoint}

For a $G$--tree $T \in \GT$ define $l_T (S) = \min_{x \in T} \max_{g
\in S} d(x,gx)$, where $S$ is some finite subset of $G$.  The
characteristic set of $S$ is $T_S = \{ x \in T \ | \ l_T(S) = \max_{g \in S}
d(x,gx) \}$.  This agrees with the earlier notion for
characteristic set when the subgroup generated by $S$ is elliptic.
Clearly for $g \in S$ we have $l_T(g) \leq l_T(S)$. We let $S'$ be the
subset of $S$ where this is an equality, i.e.\ $S' = \{g \in S \ | \
l_T(g) = l_T(S) \}$.  Finally we define $Z_S = \bigcap_{g \in S'}
T_g$.

\begin{lemma}\label{lm:T_S}
Let $T$ be a $G$--tree and let $S$ be a finite subset of $G$.  Then
$T_S$ is contained in the union of a finite simplicial tree and $Z_S$.
In particular, if $Z_S$ is a finite simplicial tree, then $T_S$ is a
finite simplicial tree.
\end{lemma}

\begin{proof}
Let $x \in T$ and $X$ be the union of all arcs from $x$ to $T_g$ for
$g \in S$, then $X$ is a finite simplicial tree.  If $y \in T_S$ is
not in $X$, let $z$ be the closest point in $X$ to $y$.  Then $d(y,gy)
\geq d(z,gz)$ for all $g \in S$ as $d(y,T_g) \geq d(z,T_g)$ with
equality only if $y \in T_g$.  If $g \in S'$ then $d(y,gy) \geq
d(z,gz) \geq l_T(S)$.  As $y \in T_S$ we have $l_T(S) \geq d(y,gy)$.
Hence we have equality $d(y,gy) = d(z,gz)$ for $g \in S'$.  Thus $y \in
T_g$ for all $g \in S'$ and hence $y \in Z_S$.
\end{proof} 

Let $S$ generate $G$.  Then for irreducible $G$--trees $T$, $Z_S$ is
finite, hence so is $T_S$.  We have some simple lemmas on the shape
and position of $T_S$ based on $l_T(S)$ and $l_T(g)$ that will be used
in Proposition \ref{prop:basepoint}.

\begin{lemma}\label{lm:T_S_shape}
Suppose that $T_S$ is finite.  Then $T_S$ is either a point or a
segment.  Moreover, the latter only occurs when there is a $g \in S$
such that $l_T(g) = l_T(S)$.  In both cases, there are distinct
$g_1,g_2 \in S$ such that $d(x,g_1 x) = d(x,g_2 x) = l_T(S)$ for all
$x \in T_S$.
\end{lemma}

\begin{proof}
Suppose $l_T(S) > \max_{g \in S} l_T(g)$ and there are distinct points
$x_1,x_2 \in T_S$.  Let $g_1,g_2 \in S$ be such that $\max_{g \in S}
d(x_i,gx_i) = d(x_i,g_i x_i) = l_T(S)$ for $i=1,2$.  Thus $x_i \notin
T_{g_i}$.  Consider the segment $[x_1,x_2]$.  Let $y \in [x_1,x_2]$
and $y \neq x_1,x_2$.  Then for any $g \in S$, $d(y,T_g) < d(x_i,T_g)$
for either $i=1$ or 2, hence $\max_{g \in S} d(y,gy) < l_T(S)$.  This
is a contradiction, therefore $T_S = \{x\}$.  Now notice that there
are $g_1,g_2 \in S$ such that $d(x,g_i x) = l_T(S)$.  For if there was
only one such $g$, then for some point $y$ near $x$ on the arc from
$x$ to $T_g$, $\max_{g \in S} d(y,gy) < d(x,gx) = l_T(S)$, which is a
contradiction.

If $l_T(S) = l_T(g)$ for $g \in S$ then $T_S \subset T_g$.  Therefore
$T_S$ is either a point or a segment.  If there were only one such $g
\in S$ such that $l_T(S) = l_T(g)$, then $T_S$ is open by a similar
argument as above.  This is a contradiction.
\end{proof}

Recall that for $A \subseteq T$ closed, we let $p_A \co T \to A$ denote
the nearest point projection.

\begin{lemma}\label{lm:T_S_pos}
Let $z \in T \setminus T_S$ and $x = p_{T_S}(z)$.  Then for some $g
\in S$ such that $d(x,gx) = l_T(S)$, we have that $x$ is on the arc
from $z$ to $T_g$.
\end{lemma}

\begin{proof}
Suppose not.  Then for points $x' \in [x,z]$ near $x$, $d(x',T_g) \leq
d(x,T_g)$ for all $g \in S$ such that $d(x,gx) = l_T(S)$.  This is a
contradiction.
\end{proof}

For an irreducible $G$--tree $T$, let $x_*$ denote the midpoint of
$T_S$.  This is called the \textit{basepoint} of the action.  Define a
map $b(T)\co G \to T$ by $g \mapsto gx_*$.  This defines a map $b : \GT
\to \GE (G,\GT )$ where $\GE (G,\GT )$ is the space of equivariant
maps from $G$ to $G$--trees in \GT.  The topology for $\GE(G,\GT)$ is
the Gromov-Hausdorff topology defined in section \ref{sc:top} where we
consider $G$ as a metric $G$-space.  The actual metric we place on $G$
does not matter as the domain is fixed in $\GE(G,\GT)$.  The remainder
of this section is used to prove that $b\co \GT \to \GE(G,\GT)$ is a
continuous function.

\begin{proposition}\label{prop:basepoint}
$b$ is continuous.
\end{proposition}

\begin{proof}
This amounts to showing that close $G$--trees in \GT \ have close
basepoints.  Let $T \in \GT$, there are two cases depending on
$l_T(S)$.

\textbf{Case 1}\qua $l_T (S) > \max_{g \in S} l_T (g)$

By Lemma \ref{lm:T_S_shape} we have that $T_S = \{ x_* \}$.  Within
the set of $g \in S$ such that $d(x_*,gx_*) = l_T(S)$, there are two
elements $g_1,g_2$ such that $x_*$ is on the spanning arc from
$T_{g_1}$ to $T_{g_2}$.  Let $x_i$ be the point on $T_{g_i}$ nearest
to $x_*$.  Thus $x_* \in [x_1,x_2]$ and $d(x_1,x_2) = d(x_1,x_*) +
d(x_*,x_2)$.

Let $U$ be the basic open set $U = U(b(T),P \times K,P,\epsilon)$,
where $S \subseteq P$ and $P( \{ x_*,x_1,x_2 \} ) \subseteq K$.  By
the remark in section \ref{sc:top}, we can assume that $P$ and $K$
contain these subsets by shrinking $U$.  Also let $V =
U(T,K,P,\delta)$, where $\delta = \frac{1}{4} \min \{ \epsilon,
d(x_*,x_1),d(x_*,x_2) \}$.  Suppose that $Y \in V$, we will show that
$b(Y) \in U$.  By definition, there is a $P$--equivariant closed
$\delta$--approximation $R \subseteq K \times L$ for some finite
subtree $L \subseteq Y$.

\begin{figure}[ht!]\small\anchor{fg:case_1}
\setlength{\unitlength}{0.0004in}
\begin{center}
{\renewcommand{\dashlinestretch}{30}
\begin{picture}(9969,2058)(0,-10)
\texture{55888888 88555555 5522a222 a2555555 55888888 88555555 552a2a2a 2a555555 
	55888888 88555555 55a222a2 22555555 55888888 88555555 552a2a2a 2a555555 
	55888888 88555555 5522a222 a2555555 55888888 88555555 552a2a2a 2a555555 
	55888888 88555555 55a222a2 22555555 55888888 88555555 552a2a2a 2a555555 }
\put(2277,1677){\shade\ellipse{108}{108}}
\put(2277,1677){\ellipse{108}{108}}
\put(8352,357){\shade\ellipse{108}{108}}
\put(8352,357){\ellipse{108}{108}}
\put(8076,957){\shade\ellipse{108}{108}}
\put(8076,957){\ellipse{108}{108}}
\put(7662,1827){\shade\ellipse{108}{108}}
\put(7662,1827){\ellipse{108}{108}}
\put(2277,942){\shade\ellipse{108}{108}}
\put(2277,942){\ellipse{108}{108}}
\put(2277,230){\shade\ellipse{108}{108}}
\put(2277,230){\ellipse{108}{108}}
\path(2277,1677)(2277,237)
\path(7662,1827)(8352,357)
\path(8277,507)(8157,252)
\path(12,1872)(14,1872)(17,1871)
	(24,1870)(35,1869)(50,1867)
	(71,1865)(97,1861)(129,1857)
	(166,1852)(210,1847)(259,1841)
	(313,1834)(371,1827)(433,1820)
	(498,1812)(566,1804)(636,1795)
	(706,1787)(778,1779)(849,1771)
	(919,1763)(988,1755)(1056,1748)
	(1123,1741)(1188,1734)(1251,1727)
	(1312,1721)(1371,1716)(1428,1711)
	(1483,1706)(1536,1701)(1588,1697)
	(1638,1694)(1687,1690)(1734,1687)
	(1781,1685)(1826,1682)(1871,1680)
	(1915,1678)(1958,1677)(2001,1676)
	(2044,1675)(2087,1675)(2087,1674)
	(2130,1674)(2173,1674)(2216,1674)
	(2259,1675)(2303,1675)(2347,1676)
	(2392,1678)(2438,1679)(2485,1681)
	(2533,1684)(2582,1686)(2633,1689)
	(2686,1693)(2740,1696)(2796,1700)
	(2854,1705)(2913,1709)(2975,1714)
	(3038,1720)(3103,1725)(3170,1731)
	(3238,1737)(3307,1744)(3376,1750)
	(3445,1757)(3514,1764)(3582,1771)
	(3648,1777)(3712,1784)(3772,1790)
	(3829,1796)(3882,1801)(3929,1806)
	(3972,1811)(4008,1815)(4039,1818)
	(4065,1821)(4085,1823)(4100,1825)
	(4110,1826)(4117,1826)(4120,1827)(4122,1827)
\path(57,57)(58,57)(62,58)
	(68,58)(78,59)(92,61)
	(110,63)(134,66)(164,69)
	(199,73)(239,77)(285,82)
	(335,87)(390,93)(449,99)
	(511,106)(576,113)(643,120)
	(712,127)(782,134)(851,141)
	(921,148)(990,155)(1058,161)
	(1125,168)(1190,174)(1254,180)
	(1316,185)(1376,191)(1434,195)
	(1490,200)(1544,204)(1597,208)
	(1647,212)(1696,215)(1744,218)
	(1790,221)(1836,224)(1880,226)
	(1923,228)(1965,230)(2007,231)
	(2048,232)(2089,233)(2129,234)
	(2169,234)(2170,234)(2213,235)
	(2257,235)(2301,235)(2345,234)
	(2390,234)(2435,233)(2481,232)
	(2527,230)(2575,228)(2623,226)
	(2673,223)(2725,220)(2778,217)
	(2833,214)(2890,210)(2948,206)
	(3009,201)(3071,196)(3135,191)
	(3200,186)(3267,180)(3334,174)
	(3402,168)(3470,162)(3537,156)
	(3604,150)(3668,144)(3729,138)
	(3787,132)(3841,127)(3890,122)
	(3934,118)(3972,114)(4005,111)
	(4032,108)(4053,106)(4069,104)
	(4080,103)(4087,103)(4090,102)(4092,102)
\path(5652,1872)(5654,1872)(5657,1871)
	(5664,1871)(5674,1869)(5689,1868)
	(5709,1866)(5735,1863)(5766,1859)
	(5803,1855)(5846,1850)(5894,1845)
	(5948,1839)(6007,1833)(6069,1826)
	(6136,1819)(6205,1811)(6276,1804)
	(6349,1796)(6423,1789)(6497,1781)
	(6571,1774)(6645,1767)(6717,1760)
	(6787,1753)(6857,1747)(6924,1741)
	(6989,1735)(7053,1730)(7114,1725)
	(7173,1721)(7230,1716)(7285,1713)
	(7338,1709)(7390,1706)(7440,1704)
	(7488,1702)(7535,1700)(7581,1698)
	(7626,1697)(7669,1696)(7712,1696)
	(7755,1696)(7797,1696)(7838,1696)
	(7879,1697)(7880,1697)(7924,1698)
	(7969,1700)(8014,1702)(8059,1704)
	(8103,1707)(8149,1710)(8195,1713)
	(8241,1717)(8289,1722)(8337,1727)
	(8387,1732)(8438,1739)(8491,1745)
	(8545,1752)(8601,1760)(8659,1768)
	(8719,1777)(8780,1787)(8842,1796)
	(8907,1807)(8972,1818)(9038,1829)
	(9104,1840)(9171,1851)(9237,1863)
	(9301,1874)(9364,1885)(9424,1896)
	(9480,1907)(9533,1916)(9580,1925)
	(9623,1933)(9661,1940)(9692,1946)
	(9718,1951)(9739,1955)(9754,1958)
	(9765,1960)(9772,1961)(9775,1962)(9777,1962)
\path(5667,12)(5669,12)(5672,13)
	(5679,14)(5689,15)(5705,17)
	(5725,19)(5751,23)(5783,27)
	(5821,32)(5865,37)(5915,43)
	(5969,50)(6029,58)(6093,66)
	(6161,74)(6232,83)(6305,92)
	(6379,101)(6455,110)(6531,119)
	(6606,128)(6681,137)(6755,146)
	(6828,154)(6898,162)(6967,170)
	(7034,178)(7099,185)(7162,192)
	(7223,198)(7281,204)(7338,210)
	(7393,215)(7446,220)(7497,225)
	(7547,229)(7595,233)(7643,237)
	(7689,241)(7734,244)(7779,247)
	(7823,250)(7866,252)(7909,255)
	(7952,257)(7997,259)(8041,261)
	(8086,263)(8130,264)(8175,265)
	(8220,266)(8266,267)(8312,268)
	(8359,268)(8408,269)(8457,269)
	(8508,268)(8560,268)(8614,267)
	(8669,267)(8726,266)(8785,264)
	(8845,263)(8907,261)(8971,260)
	(9036,258)(9102,256)(9169,254)
	(9236,251)(9303,249)(9370,246)
	(9436,244)(9500,241)(9561,239)
	(9620,237)(9675,234)(9725,232)
	(9771,230)(9812,228)(9848,227)
	(9878,226)(9902,225)(9921,224)
	(9936,223)(9946,223)(9952,222)
	(9956,222)(9957,222)
\put(402,1362){\makebox(0,0)[lb]{\smash{{{\SetFigFont{12}{24.0}{\rmdefault}{\mddefault}{\updefault}$T_{g_1}$}}}}}
\put(2457,1467){\makebox(0,0)[lb]{\smash{{{\SetFigFont{12}{24.0}{\rmdefault}{\mddefault}{\updefault}$x_1$}}}}}
\put(2457,897){\makebox(0,0)[lb]{\smash{{{\SetFigFont{12}{24.0}{\rmdefault}{\mddefault}{\updefault}$x_*$}}}}}
\put(2457,372){\makebox(0,0)[lb]{\smash{{{\SetFigFont{12}{24.0}{\rmdefault}{\mddefault}{\updefault}$x_2$}}}}}
\put(162,312){\makebox(0,0)[lb]{\smash{{{\SetFigFont{12}{24.0}{\rmdefault}{\mddefault}{\updefault}$T_{g_2}$}}}}}
\put(7797,1887){\makebox(0,0)[lb]{\smash{{{\SetFigFont{12}{24.0}{\rmdefault}{\mddefault}{\updefault}$y_1$}}}}}
\put(8262,942){\makebox(0,0)[lb]{\smash{{{\SetFigFont{12}{24.0}{\rmdefault}{\mddefault}{\updefault}$y_*$}}}}}
\put(8562,402){\makebox(0,0)[lb]{\smash{{{\SetFigFont{12}{24.0}{\rmdefault}{\mddefault}{\updefault}$y_2$}}}}}
\put(6117,1392){\makebox(0,0)[lb]{\smash{{{\SetFigFont{12}{24.0}{\rmdefault}{\mddefault}{\updefault}$Y_{g_1}$}}}}}
\put(6057,282){\makebox(0,0)[lb]{\smash{{{\SetFigFont{12}{24.0}{\rmdefault}{\mddefault}{\updefault}$Y_{g_2}$}}}}}
\end{picture}
}
\end{center}
\caption{The characteristic sets in $T$ and the
related points in $Y$ for case 1 in Proposition
\ref{prop:basepoint}}\label{fg:case_1} \end{figure}
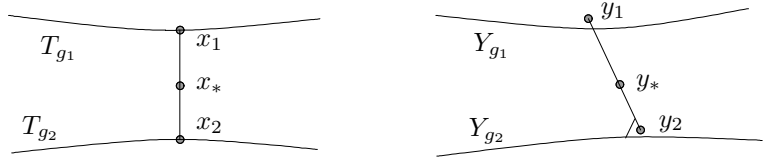

Fix related points in $L$: $x_* R y_*, x_i R y_i$ for $i=1,2$.  By
fullness of $R$, we may assume that $y_* \in [y_1,y_2]$.  Our object
now is to show that $y_*$ is close to every point in $Y_S$, in
particular, the midpoint of $Y_S$.  This involves some inequalities.
As $| d(x_*,g x_*) - d(y_*,g y_*) | < \delta$ for all $g \in S$ we
have $\max_{g \in S} d(y_*,gy_*) < \max_{g \in S} d(x_*,gx_*) + \delta
= d(x_*,g_i x_*) + \delta < d(y_*,g_i y_*) + 2 \delta$ for $i=1,2$.
Therefore, 
if $y \in Y_S$, then $d(y,g_iy) \leq \max_{g \in S} d(y,gy)
\leq \max_{g \in S} d(y_*,gy_*) < d(y_*,g_i y_*) + 2 \delta$ hence
$d(y,Y_{g_i}) < d(y_*,Y_{g_i}) + \delta$ for $i=1,2$.  We will now
show that $y,y_*$ are close to the spanning arc $\alpha$, from
$Y_{g_1}$ to $Y_{g_2}$.
$$\begin{array}{rl}\label{ar:inequal_1} | d(y_*,Y_{g_i}) - d(y_*,y_i) -
  d(y_i,Y_{g_i}) | =& |\frac{1}{2} (d(y_*,g_i y_*) - l_{Y}(g_i)) -
  d(y_*,y_i) - \\ & \frac{1}{2} (d(y_i,g_i y_i) - l_{Y}(g_i)) | \\ =&
  |\frac{1}{2}(d(y_*,g_i y_*) - d(x_*,g_i x_*)) + \\ &
  \frac{1}{2}(d(x_i,g_i x_i) - d(y_i,g_i y_i)) + \\ & (d(x_*,x_i) -
  d(y_*,y_i)) | \\ <& 2 \delta.
\end{array}$$
Hence $d(y_*,Y_{g_i}) - d(y_i,Y_{g_i}) > d(y_*,y_i) - 2\delta >
d(x_*,x_i) - 3\delta > 0$.  As $y_* \in [y_1,y_2]$ we have that $y_*$
is on $\alpha$.  Thus for $y \in Y_S$, $d(y_*,y) < \delta$.  Let $y_0
\in Y_S$ be the basepoint.

We claim that $(Id_{G},R_{\delta})$ is a $P$--equivariant closed
$\epsilon$--approximation between $b(T)\co P \to K$ and $b(Y)\co P \to
L$.  The only nontrivial check is that for $g \in P$, $b(T)(g)
R_{\delta} b(Y)(g)$.  This follows from the following calculation as
for $g \in P$ we have $gx_* R gy_*$:
 
\centerline{$d(gx_*,b(T)(g)) + d(gy_*,b(Y)(g)) = d(gy_*,gy_0) =
d(y_*,y_0) < \delta$.}

This implies $b(Y) \in U$.

\textbf{Case 2}\qua: $l_T (S) = \max_{g \in S} l_T (g)$

Let $h \in S$ be such that $l_T (h) = l_T (S)$, then $T_S \subset T_h$
as in Lemma \ref{lm:T_S_shape}.  If $x_1 \neq x_2$ assume that $h$
translates from $x_1$ to $x_2$.

Let $U$ be the basic open set $U = U(b(T),P \times K,P,\epsilon)$
where $S,S^{-1} \subseteq P$ and $P([h^{-1}x_1,hx_2]) \subseteq K$.
As in case 1, this is possible by shrinking $U$.  Let $V =
U(T,K,P,\delta)$ where $\delta = \frac{1}{9} \min \{ \epsilon, l_T
(S) \}$.  Suppose that $Y \in V$, we will show that $b(Y) \in U$.  By
definition, there is a $P$--equivariant closed $\delta$--approximation
$R \subseteq K \times L$ for some finite subtree $L \subseteq Y$.

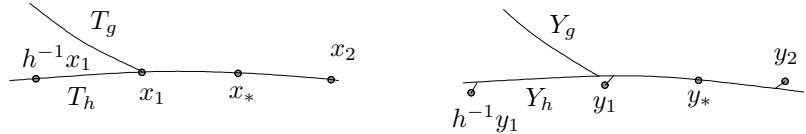
\begin{figure}[ht!]\small\anchor{fg:case_2}
\setlength{\unitlength}{0.0004in}
\begin{center}
{\renewcommand{\dashlinestretch}{30}
\begin{picture}(10494,1721)(0,-10)
\texture{55888888 88555555 5522a222 a2555555 55888888 88555555 552a2a2a 2a555555 
	55888888 88555555 55a222a2 22555555 55888888 88555555 552a2a2a 2a555555 
	55888888 88555555 5522a222 a2555555 55888888 88555555 552a2a2a 2a555555 
	55888888 88555555 55a222a2 22555555 55888888 88555555 552a2a2a 2a555555 }
\put(357,695){\shade\ellipse{84}{84}}
\put(357,695){\ellipse{84}{84}}
\put(3025,765){\shade\ellipse{84}{84}}
\put(3025,765){\ellipse{84}{84}}
\put(10242,660){\shade\ellipse{84}{84}}
\put(10242,660){\ellipse{84}{84}}
\put(9102,675){\shade\ellipse{84}{84}}
\put(9102,675){\ellipse{84}{84}}
\put(1755,780){\shade\ellipse{84}{84}}
\put(1755,780){\ellipse{84}{84}}
\put(4250,685){\shade\ellipse{84}{84}}
\put(4250,685){\ellipse{84}{84}}
\put(6100,510){\shade\ellipse{84}{84}}
\put(6100,510){\ellipse{84}{84}}
\put(7865,610){\shade\ellipse{84}{84}}
\put(7865,610){\ellipse{84}{84}}
\path(10107,570)(10242,660)
\path(6177,645)(6102,510)
\path(7977,735)(7872,615)
\path(12,675)(13,675)(17,675)
	(23,676)(33,677)(47,678)
	(66,679)(90,681)(120,683)
	(155,686)(196,689)(242,692)
	(293,696)(349,700)(409,705)
	(472,709)(538,714)(606,719)
	(676,724)(747,729)(818,734)
	(889,739)(959,744)(1028,748)
	(1097,753)(1163,757)(1228,761)
	(1292,765)(1353,768)(1413,772)
	(1470,775)(1526,778)(1580,781)
	(1633,783)(1684,785)(1733,788)
	(1782,789)(1829,791)(1875,792)
	(1920,794)(1965,795)(2009,796)
	(2052,796)(2096,797)(2139,797)
	(2182,797)(2229,798)(2276,797)
	(2324,797)(2371,797)(2419,796)
	(2467,795)(2516,795)(2565,794)
	(2614,792)(2663,791)(2713,789)
	(2763,788)(2813,786)(2863,784)
	(2913,782)(2963,780)(3013,778)
	(3063,775)(3112,773)(3161,771)
	(3210,768)(3257,765)(3304,763)
	(3350,760)(3396,757)(3440,755)
	(3483,752)(3525,749)(3566,746)
	(3605,744)(3643,741)(3680,738)
	(3716,736)(3750,733)(3783,731)
	(3814,728)(3845,726)(3874,723)
	(3901,721)(3928,719)(3953,717)
	(3977,715)(4029,711)(4075,707)
	(4116,703)(4152,700)(4184,697)
	(4212,694)(4237,691)(4259,688)
	(4279,686)(4297,683)(4312,681)
	(4324,679)(4333,678)(4340,676)
	(4344,676)(4346,675)(4347,675)
\path(1763,794)(1761,795)(1755,797)
	(1746,802)(1731,809)(1711,818)
	(1685,830)(1654,845)(1619,862)
	(1580,880)(1538,900)(1495,920)
	(1451,941)(1408,962)(1366,982)
	(1325,1002)(1287,1021)(1250,1039)
	(1215,1057)(1182,1074)(1151,1090)
	(1121,1106)(1093,1121)(1066,1136)
	(1040,1150)(1015,1165)(990,1179)
	(966,1194)(941,1209)(917,1224)
	(893,1239)(868,1255)(843,1272)
	(818,1289)(791,1308)(764,1327)
	(735,1348)(705,1370)(673,1393)
	(640,1417)(606,1443)(570,1470)
	(534,1497)(498,1524)(462,1552)
	(428,1578)(396,1603)(367,1625)
	(341,1645)(320,1661)(304,1674)
	(292,1683)(284,1689)(280,1693)(278,1694)
\path(5997,645)(5998,645)(6002,645)
	(6008,646)(6017,646)(6031,647)
	(6050,648)(6074,650)(6103,652)
	(6138,654)(6179,656)(6225,659)
	(6276,662)(6333,665)(6393,669)
	(6458,672)(6525,676)(6596,680)
	(6668,684)(6742,688)(6816,692)
	(6891,696)(6966,700)(7040,704)
	(7112,708)(7184,711)(7254,714)
	(7322,718)(7389,720)(7454,723)
	(7516,726)(7577,728)(7635,730)
	(7692,732)(7747,733)(7801,735)
	(7852,736)(7903,737)(7952,738)
	(7999,738)(8046,738)(8092,739)
	(8137,738)(8181,738)(8225,738)
	(8268,737)(8312,736)(8354,735)
	(8355,735)(8399,734)(8444,732)
	(8488,731)(8533,729)(8578,726)
	(8623,724)(8669,721)(8716,718)
	(8763,715)(8811,712)(8860,708)
	(8911,704)(8963,699)(9017,694)
	(9072,689)(9128,684)(9187,678)
	(9247,672)(9308,665)(9372,658)
	(9436,651)(9503,644)(9570,636)
	(9638,628)(9707,620)(9775,612)
	(9844,604)(9911,596)(9976,588)
	(10040,580)(10101,573)(10158,566)
	(10212,559)(10261,553)(10305,548)
	(10345,543)(10378,538)(10407,535)
	(10430,532)(10448,529)(10462,528)
	(10472,526)(10477,526)(10481,525)(10482,525)
\path(7787,729)(7785,730)(7780,734)
	(7770,739)(7756,748)(7736,760)
	(7710,776)(7680,794)(7645,815)
	(7606,839)(7566,864)(7524,889)
	(7481,915)(7439,941)(7398,966)
	(7359,991)(7322,1014)(7287,1036)
	(7253,1057)(7222,1076)(7193,1095)
	(7166,1112)(7141,1129)(7116,1145)
	(7093,1160)(7071,1175)(7050,1190)
	(7030,1204)(7004,1222)(6979,1240)
	(6954,1258)(6930,1276)(6906,1294)
	(6882,1312)(6859,1330)(6836,1348)
	(6813,1366)(6792,1383)(6771,1400)
	(6751,1417)(6732,1433)(6713,1448)
	(6697,1462)(6681,1476)(6666,1489)
	(6653,1500)(6641,1511)(6629,1521)
	(6619,1531)(6610,1539)(6591,1555)
	(6577,1569)(6564,1580)(6554,1589)
	(6545,1597)(6538,1604)(6532,1609)
	(6529,1612)(6527,1614)
\put(747,345){\makebox(0,0)[lb]{\smash{{{\SetFigFont{12}{24.0}{\rmdefault}{\mddefault}{\updefault}$T_h$}}}}}
\put(1047,1350){\makebox(0,0)[lb]{\smash{{{\SetFigFont{12}{24.0}{\rmdefault}{\mddefault}{\updefault}$T_g$}}}}}
\put(147,855){\makebox(0,0)[lb]{\smash{{{\SetFigFont{12}{24.0}{\rmdefault}{\mddefault}{\updefault}$h^{-1}x_1$}}}}}
\put(1692,420){\makebox(0,0)[lb]{\smash{{{\SetFigFont{12}{24.0}{\rmdefault}{\mddefault}{\updefault}$x_1$}}}}}
\put(2877,435){\makebox(0,0)[lb]{\smash{{{\SetFigFont{12}{24.0}{\rmdefault}{\mddefault}{\updefault}$x_*$}}}}}
\put(4212,990){\makebox(0,0)[lb]{\smash{{{\SetFigFont{12}{24.0}{\rmdefault}{\mddefault}{\updefault}$x_2$}}}}}
\put(6792,330){\makebox(0,0)[lb]{\smash{{{\SetFigFont{12}{24.0}{\rmdefault}{\mddefault}{\updefault}$Y_h$}}}}}
\put(7107,1305){\makebox(0,0)[lb]{\smash{{{\SetFigFont{12}{24.0}{\rmdefault}{\mddefault}{\updefault}$Y_g$}}}}}
\put(5817,60){\makebox(0,0)[lb]{\smash{{{\SetFigFont{12}{24.0}{\rmdefault}{\mddefault}{\updefault}$h^{-1}y_1$}}}}}
\put(7677,315){\makebox(0,0)[lb]{\smash{{{\SetFigFont{12}{24.0}{\rmdefault}{\mddefault}{\updefault}$y_1$}}}}}
\put(10062,945){\makebox(0,0)[lb]{\smash{{{\SetFigFont{12}{24.0}{\rmdefault}{\mddefault}{\updefault}$y_2$}}}}}
\put(8922,390){\makebox(0,0)[lb]{\smash{{{\SetFigFont{12}{24.0}{\rmdefault}{\mddefault}{\updefault}$y_*$}}}}}
\end{picture}
}\end{center}
\vspace{-.5cm}
\caption{The characteristic sets in $T$ and the
related points in $Y$ for case 2 in Proposition
\ref{prop:basepoint}}\label{fg:case_2} \end{figure}

Fix related points in $L$: $x_* R y_*, x_i R y_i$ for $i=1,2$.  Again,
by the fullness of $R$ we may assume that $y_* \in [y_1,y_2]$.  Our
object now is to show that the Hausdorff distance between $[y_1,y_2]$
and $Y_S$ is small.  As before, this involves some inequalities.  Our
first step is to show that points in $[y_1,y_2]$ are close to some
point in $Y_S$.

Let $z \in [y_1,y_2]$ and $x \in [x_1,x_2] = T_S$ where
$xRz$.  Then $\max_{g \in S} d(z,gz) < \max_{g \in S} d(x,gx) + \delta
= d(x,hx) + \delta < d(z,hz) + 2 \delta$.  Since $R$ is full, this is
true for any $z \in [y_1,y_2]$.  Note that the above inequality implies
$l_Y(S) < l_T(S) + \delta$.

We now show that the segment $[h^{-1}y_1,hy_2]$ is close to the axis
$Y_h$.  Let $z,z' \in [h^{-1}y_1,hy_2]$ and $x,x' \in [h^{-1}x_1,hx_2]$
where $xRz,x'Rz'$.
$$\begin{array}{rl}\label{ar:inequal_2} | d(z,Y_h) - d(z',Y_h) | =&
  |\frac{1}{2} (d(z,hz) - l_{Y}(h)) - \frac{1}{2}(d(z',hz') -
  l_{Y}(h)) | \\ =& \frac{1}{2}|(d(z,hz) - d(x,hx)) + (d(x',hx') -
  d(z,hz'))| \\ <& \delta.
\end{array}$$
In particular $|d(h^{-1}y_1,Y_h) - d(hy_2,Y_h) | < \delta$, as
$d(h^{-1}y_1,hy_2) > 2 l_T (h) - \delta > 2 \delta$ this implies that
there is a $z_0 \in [h^{-1}y_1,hy_2] \cap Y_h$.  Hence for any $z \in
[h^{-1}y_1,hy_2]$ we have $d(z,Y_h) < \delta$.  Likewise the same is
true for $z \in [y_1,y_2]$.

Now for $z \in [y_1,y_2]$, $l_Y(S) - 2 \delta \leq \max_{g \in S}
d(z,gz) - 2 \delta < d(z,hz) < l_Y (h) + 2 \delta < l_Y (S) + 2
\delta$.  For $z \in [y_1,y_2]$ that are not in $Y_S$, let $y =
p_{Y_S}(z)$ and let $g' \in S$ be given by Lemma \ref{lm:T_S_pos}.
Then $d(z,y) = \frac{1}{2}(d(z, g'z) - d(y,g' y)) \leq
\frac{1}{2}(\max_{g \in S} d(z,gz) - l_Y (S)) < \delta$.  Hence for
$z \in [y_1,y_2]$, we have $d(z,Y_S) < \delta$.

For the opposite inequality we show that points in $Y_S$ are close to
some point in $[y_1,y_2]$.  We do so by showing that points far enough
away from $[y_1,y_2]$ cannot lie in $Y_S$.  First note that the above
inequality implies: $l_Y(S) - l_Y (h) < 4 \delta$.  Hence if $y' \in
Y_S$, then $2d(y',Y_h) = d(y',hy') - l_Y (h) < l_Y (S) - (l_Y(S) - 4
\delta)$.  Thus $d(y',Y_h) < 2 \delta$.  Recall that we have shown
$l_Y(S) < l_T(S) + \delta$.

The idea now is to use Lemma \ref{lm:T_S_pos} on points far from
$[y_1,y_2]$.  Assume that $y' \in Y_S$ and $d(y',[y_1,y_2]) \geq
4\delta$.  Then there is some point $y \in Y_h \cap L$ with
$d(y,[y_1,y_2]) \geq 2\delta$.  Without loss of generality, we assume
that $y$ is closer to $y_1$ than to $y_2$.  Let $x \in T_h \cap K$ be
such that $xRy$.  Then $d(x,x_1) \geq \delta$.  Hence by Lemma
\ref{lm:T_S_pos} there is a $g \in S$ such that $d(x,gx) \geq l_T(S) +
2\delta$.  Therefore $l_Y(S) \geq d(y,gy) \geq l_T(S) + \delta >
l_Y(S)$, which is a contradiction.  Therefore the Hausdorff distance
between $Y_S$ and $[y_1,y_2]$ is less than $4 \delta$.  Let $y_0 \in
Y_S$ be the basepoint, then $d(y_0,y_*) < 4 \delta$.  Now proceed as
in case 1 using the $P$--equivariant closed $\epsilon$--approximation
$(Id_G,R_{4\delta})$.

This completes the proof.
\end{proof}

\begin{remark}\label{rm:basepoint}{\rm
The technical statement proved in the above which is used later on in
Lemma \ref{lm:step1} is that if two trees $Y$ and $Z$ have subtrees,
$L \subseteq Y$, $M \subseteq Z$ with $P\{b(Y)(1)\} \subseteq L$, $S
\subseteq P$ and a $P$--equivariant $\epsilon$--approximation $R
\subseteq L \times M$, then if $z \in Z$ with $b(Y)(1)Rz$, we have
$d(z,b(Z)(1))< 4 \epsilon$.  In other words, any point related to the
basepoint of $Y$ is within $4\epsilon$ of the basepoint of $Z$.  }
\end{remark}

%%%%%%%%%%%%%%%%%%%%%%%%%%%%%%%%%%%%%%%%%%%%%%%%%%%%%%%%%%%%%%%%%%%%%%%%%%%%% 

\section{Contractibility of deformation space}\label{sc:cont_def}

To prove the contractibility of the unnormalized deformation space
\GT, we construct a homotopy onto a contractible subset.  To define
the homotopy, for any $G$--tree $T' \in \GT$ we need to build a nice
map from some fixed $G$--tree $T \in \GT$ to $T'$.  To ensure that the
map $T \to T'$ is nice, we will need $T$ to be reduced.

\begin{definition}\label{def:reduced} {\rm
A $G$--tree $T$ is \textit{reduced} if for all edges $e=[u,v]$,
$u$ is $G$--equivalent to $v$ if $G_e = G_u$.}
\end{definition}

This is equivalent to Forester's definition in \cite{ar:F1} where a
tree is said to be reduced if it admits no collapse moves.  We will
use this notion via the next lemma.

\begin{lemma}\label{lm:G_equiv}
Let $T$ be a reduced $G$--tree and $u,v \in T$ vertices such that there
is an edge $e=[u,v]$ and $x$ a vertex with $G_u,G_v \subseteq G_x$.
Then $u$ is $G$--equivalent to $v$.
\end{lemma}

\begin{proof}
Without loss of generality, assume that $v$ is closer to $x$ than $u$
is.  Then $[u,x] = e \cup [v,x]$ and as $G_u$ stabilizes $[u,x]$ this
implies that $G_u = G_e$.  Hence as $T$ is reduced, the two endpoints
of $e$ are $G$--equivalent.
\end{proof}

We now require that our unnormalized deformation space \GT \ contains
a $G$--tree with finitely generated vertex groups.  In particular as
all $G$--trees in \GT \ are cocompact, there is a reduced tree $T \in
\GT$ with finitely generated vertex groups.  Define $\GS (T, \GT)$ as
the space of all continuous maps from $T$ to $G$--trees in \GT \ that
take vertices to vertices and are injective on the edges of $T$.  We
call such maps \textit{transverse}.  This has a different meaning than
in \cite{ar:S}, where transverse only implies cellular.  We topologize
$\GS(T,\GT)$ using the Gromov-Hausdorff topology from section
\ref{sc:top}.

Our aim now is to build a section $B\co \GT \to \GS (T,\GT)$.  Let $G$
be finitely generated by $S$ and fix $X \subseteq T$ a subtree whose
edges map bijectively to $T/G$.  We follow Forester's construction
from Proposition 4.16 in \cite{ar:F1}.  Order the vertices of $X$ as
$\{ v_1, \ldots ,v_k \}$ where vertices in the same orbit are
consecutive.  For the $i$th orbit $v_{i_0}, \ldots, v_{i_0 + d}$ let
$g_{i_0} = 1$ and fix $g_{i_0 + q} \in G$ such that $g_{i_0 +
q}v_{i_0} = v_{i_0 +q}$ for $1 \leq q \leq d$.  As the path $[v_{i_0},
g^{-1}_{i_0 + q}g_{i_0 + p}v_{i_0}]$ for $1\leq p,q \leq d$ is
contained in $X$, it maps bijectively to $T/G$.  Therefore the
products $g^{-1}_{i_0 + q}g_{i_0 + p}$ are hyperbolic for $p \neq q$
(Lemma 2.7(b) \cite{ar:F1}).

Given $Y \in \GT$, we define the map $B(Y)\co T \to Y$ first on the
vertices of $X$. Let $y_* = b(Y)(1)$, where $b$ is the basepoint map
of Proposition \ref{prop:basepoint}.  Recall that $p_A$ is projection
onto the closed subset $A$.  Consider an orbit
$\{v_{i_0},\ldots,v_{i_0 + d} \}$.  If $G_{v_{i_0}} \neq \{1\}$ then
let $Y_{i_0} \subseteq Y$ be the characteristic set of $G_{v_{i_0}}$.
Otherwise, let $Y_{i_0} = y_*$.  As $T$ is reduced, $G_{v_{i_0}} =
\{1\}$ can only happen if $G$ is a finitely generated free group of
rank at least 2.  In which case $T/G$ is a rose and there is only one
orbit of vertices in $X$.  Define $B(Y)$ on the orbit by: $v_{i_0 + d}
\mapsto g_{i_0 + d}p_{Y_{i_0}}(y_*)$.

We now show that $B(Y)$ can be extended to a transverse map.  If there
is an edge $e \subseteq X$ where $e = [u,v]$ with $B(Y)(u) = B(Y)(v) =
x' \in Y $, then $G_u,G_v \subseteq G_{x'}$.  This subgroup must fix a
vertex $x \in T$, hence $G_u,G_v \subseteq G_x$ and by Lemma
\ref{lm:G_equiv}, $u$ and $v$ must be in the same orbit.  But if $v_i$
and $v_j$ are in the same orbit then as ${g_i}^{-1}g_j$ is hyperbolic
for $i \neq j$ necessarily $B(Y)(v_i) \neq B(Y)(v_j)$.  Thus we can
linearly map each edge of $X$ injectively into $Y$.  Now extend $B(Y)$
to all of $T$ equivariantly.  As $B(Y)$ is injective on each edge this
defines $B\co \GT \to \GS (T,\GT)$.

For the $i$th orbit, let $G_i$ be the vertex stabilizer of the first
vertex in this orbit and denote the characteristic set for $G_i$ by
the subscript $i$, i.e.\ $Y_{G_i} = Y_i$.  If $G_i = 1$ then as before,
set $Y_i = y_*$.  Let $G_i$ be finitely generated by $S_i$, then for
any $G$--tree $Y \in \GT$ we have $Y_i = \cap_{s \in S_i} Y_s$.  Let $Q
\subseteq G$ be the union of the $S_i$'s and $S$, a finite generating
set for $G$.

\begin{lemma}\label{lm:range-section}
$B$ is continuous and $\Ra(B(Y)) = Y$ for all $Y \in \GT$.
\end{lemma}

If $G$ is finitely generated free group of rank at least 2, then this
follows from Proposition \ref{prop:basepoint}.  Thus we assume that
$G$ is not free.  Before we prove this lemma in general, we prove a
statement about the position of the basepoint relative the fixed point
sets.

\begin{lemma}\label{lm:step1}
Let $Y,Z \in \GT$, and $Y_i,Z_i$ be the characteristic sets as
described above.  Let $L \subseteq Y,M \subseteq Z$ be subtrees with
$R \subseteq L \times M$ a $P$--equivariant $\epsilon$--approximation
where $Q \subseteq P$ and $P \{y_*,p_{Y_i}(y_*) \} \subseteq L$.  Then
for $z \in Z$ such that $p_{Y_i}(y_*) R z$ we have $d(z,p_{Z_i}(z_*))
< 10 \epsilon$.
\end{lemma}

\begin{proof}
Fix $\hat{z} \in Z$ where $y_* R \hat{z}$, then by Remark
\ref{rm:basepoint}, $d(\hat{z},z_*) < 4 \epsilon$.  The lemma follows
from 3 observations:

\noindent
(i)\qua  $|d(z_*,z) - d(y_*,p_{Y_i}(y_*)| < 
|d(\hat{z},z) - d(y_*,p_{Y_i}(y_*))| + 4 \epsilon < 5 \epsilon$.

\noindent
(ii)\qua Let $g \in S_i$ be such that $p_{Z_i}(z_*) = p_{Z_g}(z_*)$.  Then:
$$\begin{array}{rl} 2d(z_*,p_{Z_i}(z_*)) = d(z_*,gz_*) <&
d(\hat{z},g\hat{z}) + 8 \epsilon \\ <& d(y_*,gy_*) + 9\epsilon =
2d(y_*,p_{Y_g}(y_*)) + 9\epsilon \\ \leq& 2d(y_*,p_{Y_i}(y_*)) +
9\epsilon.
\end{array}$$
Likewise, running this argument with $h \in S_i$ such that
$p_{Y_i}(y_*) = p_{Y_h}(y_*)$, we see that $|d(z_*,p_{Z_i}(z_*)) -
d(y_*,p_{Y_i}(y_*))| < 5\epsilon$.

\noindent
(iii)\qua Let $g \in S_i$ be such that $p_{Z_i}(z) = p_{Z_g}(z)$. Then:
$2d(z,p_{Z_i}(z)) = d(z,gz) < d(p_{Y_i}(y_*),gp_{Y_i}(y_*)) + 
\epsilon = \epsilon$.

Putting (i) and (ii) together: $|d(z_*,z) - d(z_*,p_{Z_i}(z_*))| < 
10\epsilon$.  If $[z_*,z]$ passes through $Z_i$ then 
$d(z,p_{Z_i}(z_*)) < 10\epsilon$.  If $[z_*,z]$ doesn't pass through 
$Z_i$, then $p_{Z_i}(z) = p_{Z_i}(z_*)$ and hence $d(z,p_{Z_i}(z_*)) =
d(z,p_{Z_i}(z))< 
\epsilon$ by (iii).
\end{proof}

Now we can prove Lemma \ref{lm:range-section}.

\begin{proof}
Let $Y \in \GT$ and let $U$ be a basic open set of $B(Y),U=U(B(Y),K
\times L,P, \epsilon)$ where $Q \subseteq P$ and $P \{y_*,p_{Y_i}(y_*)
\} \subseteq L$.  Enlarge $P$ such that $K \subseteq PX$.  Also let
$V$ be a basic open set for $Y \in \GT$, $V=U(Y,L,P,\delta)$ where
$\delta = \frac{\epsilon}{21}$.  If $Z \in V$, we have by definition a
$\delta$--approximation $R \subseteq L \times M$ for some $M \subseteq
Z$.  This is the set-up in Lemma \ref{lm:step1}.  As before, we will
show that $B(Z) \in U$.

We claim that $(Id,R_{10\delta})$ is an $\epsilon$--approximation from
$B(Y)$ to $B(Z)$, hence $B(Z) \in U$.  As in the proof of Proposition
\ref{prop:basepoint}, the only nontrivial check is that $B(Y)(x)
R_{10\delta} B(Z)(x)$ for $x \in K$.  Without loss of generality, we
can assume that $x$ is a vertex as the maps are linear on the edges.
Now $x=gv$ for some $v \in X$, ordered first in its orbit and some $g
\in P$.  Let $A = Y_{G_v},B=Z_{G_v}$, then $B(Y)(x) = gp_A(y_*)$ and
$B(Z)(x) = gp_B(z_*)$.  Let $z \in Z$ be such that $p_A(y_*) R z$.
Thus by Lemma \ref{lm:step1}:

\centerline{$d(gp_A(y_*),B(Y)(x)) + d(gz,B(Z)(x))= d(gz,gp_B(z_*)) =
d(z,p_B(z_*)) < 10 \delta$.}

\noindent
This completes the proof.
\end{proof}

The map $B(Y)$ is not a morphism in the sense used within this paper.
However we can redefine the metric on $T$ to get a new $G$--tree $T_Y$
such that the map $B(Y)$ when regarded as a map $B(Y)\co T_Y \to Y$ is a
morphism.  As each edge of $T$ is mapped injectively via $B(Y)$ to $Y$
we can remetrize each edge by pulling back the metric on $Y$.  Thus we
can remetrize $T$ by setting the distance between two points to be the
length of the geodesic path between them.  Call this new $G$--tree
$T_Y$.  Then $T_Y$ is equivariantly homeomorphic to the $G$--tree $T$.
Let $\GL(T)$ denote the set of $G$--trees in $\GT$ which are equivariantly
homeomorphic to the $G$--tree $T$.  Recall that the volume of $T$ is
defined as $\vol(T) = \sum \length(e)$ where the sum is over the
unoriented edges of $T/G$.

\begin{proposition}\label{prop:T_Y}
$\GL(T)$ is homeomorphic to $\sigma \times \R^+$ where $\sigma$ is an open
simplex of dimension one less than the number of edges of $T/G$.
\end{proposition}

\begin{proof}
Fix an ordering $e_1,\ldots, e_n$ of the edges of $T$.  This in turn
gives an ordering of the edges of $T' \in \GL(T)$.  Let $h: \GL(T) \to
\sigma \times \R^+$ be defined by:
\begin{equation}\label{eq:homeo}
h(T') =
(\frac{1}{\vol(T')}(\length(e_1), \ldots,\length(e_n)),\vol(T')). 
\end{equation}
It is clear that this map gives a bijection between the sets.  As we
are working with irreducible $G$--trees, as mentioned in section
\ref{sc:top} the Gromov--Hausdorff topology is the same as the axes
topology.  A small change in $\GL(T)$ of the length functions results in
a small change in the lengths of the edges.  And conversely, a small
change in the length of the edges of $T/G$ results in a small change
of the length functions only for the hyperbolic conjugacy classes
whose axis project down to paths which cross the rescaled edges.
Therefore, $h$ is a homeomorphism.
\end{proof}

Denote by $\beta(Y)\co T_Y \to Y$ the morphism induced by the
transverse map $B(Y)\co T \to Y$.  Hence $\beta$ defines a map
$\beta\co \GT \to \GM(\GL(T),\GT)$.  As $B\co \GT \to \GS(T,\GT)$ is
continuous and the newly defined metric on $T_Y$ depends continuously
on the metric on $Y$, $\beta$ is continuous.  We can now define a
homotopy equivalence from \GT \ to $\GL(T)$.

\begin{theorem}\label{th:X_contractible}
For a finitely generated group $G$, any irreducible unnormalized
deformation space which contains a $G$--tree with finitely generated
vertex groups is contractible.
\end{theorem}

\begin{proof}
Let $\beta\co \GT \to \GM(\GL(T),\GT),\Phi\co \GM (\GT) \times \{
(s,t) \ | 0 \leq s \leq t \leq 1 \} \to \GM (\GT)$ and $\Ra\co
\GM(\GT) \to \GT$ be the continuous functions defined above.

Define a homotopy $H\co \GT \times [0,1] \to \GT$ by $H_{(1-t)}(Y) =
\Ra(\Phi(\beta(Y),0,t))$.  Then $H_0(Y) = \Ra(\Phi(\beta(Y)),0,1)) =
\Ra(\beta(Y)) = Y$ and $H_1(\GT) = \GL(T)$, which is contractible by
\ref{prop:T_Y}.
\end{proof}

Recall that $\GD = \GT / \R^{+}$.  As $\Ra \circ \Phi$ is
$\R^+$--invariant (Remark \ref{rm:R-invar}) and $\beta$ clearly is
also, $H$ descends to a homotopy of \GD.  Therefore we have the
following theorem as stated in the introduction:

\begin{theorem}\label{th:d_contractible}
For a finitely generated group $G$, any irreducible deformation space
which contains a $G$--tree with finitely generated vertex groups is
contractible.
\end{theorem}

%%%%%%%%%%%%%%%%%%%%%%%%%%%%%%%%%%%%%%%%%%%%%%%%%%%%%%%%%%%%%%%%%%%%%%%%%%%%% 

%%%%%%%%%%%%%%%%%%%%%%%%%%%%%%%%%%%%%%%%%%%%%%%%%%%%%%%%%%%%%%%%%%%%%%%%%%%%% 

\Addressesr


\begin{thebibliography}{99}

\bibitem{ar:B} \textbf{M Bestvina}, \textit{The topology of $\Out(F_n)$},
Proc. of the ICM, Beijing II, Higher Ed. Press, Beijing (2002) 373--384
\MR{1957048}


% MR lookup returned no matches
\bibitem{ar:BV} \textbf{M Bridson}, \textbf{K Vogtmann},
\textit{Automorphism groups of free groups, surface groups and free
abelian groups}, 
%to appear in ``Problems on mapping class groups and related topics'', 
%Proc. Symp. Pure and Applied Math. (B Farb, editor),
\arxiv{math.GR/0507612}

\bibitem{ar:CM}
\textbf{M Culler}, \textbf{J\,W Morgan}, \emph{Group actions on {${\bf
  R}$}-trees}, Proc. London Math. Soc. (3) 55 (1987) 571--604 \MR{907233}

\bibitem{ar:CV}
\textbf{M Culler}, \textbf{K Vogtmann}, \emph{Moduli of graphs and
  automorphisms of free groups}, Invent. Math. 84 (1986) 91--119 \MR{830040}

\bibitem{ar:F1}
\textbf{M Forester}, \emph{Deformation and rigidity of simplicial group actions
  on trees}, \gtref6{2002}{8}{219}{267} \MR{1914569}

\bibitem{ar:G}
\textbf{M Gromov}, \emph{Groups of polynomial growth and expanding maps}, Inst.
  Hautes \'Etudes Sci. Publ. Math.  (1981) 53--73 \MR{623534}

\bibitem{ar:GL} \textbf{V Guirardel}, \textbf{G Levitt}, \textit{A general
construction of JSJ splittings}, research announcement

\bibitem{ar:GL2} \textbf{V Guirardel}, \textbf{G Levitt}, \textit{The
outer space of a free product}, e-print (2005)
\arxiv{math.GR/0501288}

\bibitem{ar:MM}
\textbf{D McCullough}, \textbf{A Miller}, \emph{Symmetric automorphisms of
  free products}, Mem. Amer. Math. Soc. 122 (1996) viii+97pp \MR{1329943}

\bibitem{ar:P}
\textbf{F Paulin}, \emph{The {G}romov topology on {${\bf R}$}-trees}, Topology
  Appl. 32 (1989) 197--221 \MR{1007101}

\bibitem{ar:S} \textbf{R Skora}, \textit{Deformation of length functions in
groups}, preprint

\bibitem{ar:V}
\textbf{K Vogtmann}, \emph{Automorphisms of free groups and outer space}, Geom.
  Dedicata 94 (2002) 1--31 \MR{1950871}

\end{thebibliography}
\end{document}